\documentclass[%
 aip,
 amsmath,amssymb,
 reprint,%
]{revtex4-1}
\usepackage{graphicx}
\usepackage{dcolumn}
\usepackage{bm}
\usepackage[utf8]{inputenc}
\usepackage[T1]{fontenc}
\usepackage{mathptmx}
\usepackage{etoolbox}
\usepackage{bbm}
\usepackage{xcolor}

\renewcommand{\vec}[1]{\mathbf{#1}}
\renewcommand{\d}{{\rm d}}

\newcommand{\PD}[2]{\frac{\partial #1}{\partial #2}}
\newcommand{\FD}[2]{\frac{\d #1}{\d #2}}

\draft 

\makeatletter
\def\@email#1#2{%
 \endgroup
 \patchcmd{\titleblock@produce}
  {\frontmatter@RRAPformat}
  {\frontmatter@RRAPformat{\produce@RRAP{*#1\href{mailto:#2}{#2}}}\frontmatter@RRAPformat}
  {}{}
}%
\makeatother

\begin{document}

\title{Insights into oscillator network dynamics using a phase-isostable framework} 

\author{R. Nicks}
\email{Rachel.Nicks@nottingham.ac.uk}
\author{R. Allen}
\author{S. Coombes}

\affiliation{School of Mathematical Sciences, University of Nottingham, Nottingham, NG7 2RD, UK}

\date{\today}

\begin{abstract}
Networks of coupled nonlinear oscillators can display a wide range of emergent behaviours under variation of the strength of the coupling. Network equations for pairs of coupled oscillators where the dynamics of each node is described by the evolution of its phase and slowest decaying isostable coordinate have previously been shown to capture bifurcations and dynamics of the network which cannot be explained through standard phase reduction. An alternative framework using isostable coordinates to obtain higher-order phase reductions has also demonstrated a similar descriptive ability for two oscillators. In this work we consider the phase-isostable network equations for an arbitrary but finite number of identical coupled oscillators, obtaining conditions required for stability of phase-locked states including synchrony. For the mean-field complex Ginzburg-Landau equation where the solutions of the full system are known, we compare the accuracy of the phase-isostable network equations and higher-order phase reductions in capturing bifurcations of phase-locked states.  We find the former to be the more accurate and therefore employ this to investigate the dynamics of globally linearly coupled networks of Morris-Lecar neuron models  (both two and many nodes). We observe qualitative correspondence between results from numerical simulations of the full system and the phase-isostable description demonstrating that in both small and large networks the phase-isostable framework is able to capture dynamics that the first-order phase description cannot. 
\end{abstract}

\maketitle 

\begin{quotation}
The utility of the classical technique of phase reduction for describing the dynamics of networks of coupled oscillators is limited by the assumption that the dynamics for each node remain on the stable limit cycle of the uncoupled system. We here investigate reduced equations for networks of arbitrary finite size where the dynamics of each node is described in terms of its phase and the slowest decaying isostable coordinate, allowing for representation of trajectories away from (but near) the limit cycle. Specifically, we consider conditions for the existence and stability of phase-locked states generalising existing results for phase-reduced equations. We show that phase-isostable network equations provide the most accurate qualitative description of dynamics of the mean-field Ginzburg-Landau equation when compared with alternative frameworks such as higher-order phase reduction \cite{Park2021, Leon2019}. We further demonstrate the power of the general framework by considering networks of neural oscillators. We observe phenomena including the emergence of quasiperiodic behaviour that cannot be captured using first-order phase reduction. The results are shown to be in good qualitative agreement with the  dynamics of the original network through numerical simulations and bifurcation analysis. 
\end{quotation}

\section{Introduction}

Oscillations observed in biological, physical and chemical systems are often due to the presence of attracting limit cycles within the high dimensional dynamics \cite{Kuramoto1984, Winfree2001}. The classical technique of first-order phase reduction \cite{Ermentrout2010, Brown2004, Izhikevich2007, Pietras2019, Nakao2016, Monga2019} provides a rigorous way of describing the dynamics of weakly perturbed oscillators in terms of a single phase variable using the notion of isochrons that extend the phase variable for a limit-cycle attractor to its basin of attraction \cite{Winfree1967, Guckenheimer1975}. The power of this approach has been demonstrated through its ability to reveal complex dynamics of weakly forced oscillations and emergent behaviours in weakly coupled oscillator networks in a variety of relevant systems \cite{Dorfler2013, Dorfler2014, Abrams2004, Hoppensteadt1997, Pietras2019}. However, phase reduction assumes that the dynamics remain close to the unperturbed limit cycle and therefore requires that interactions are weak and convergence to the limit cycle is fast. The ability of phase reduced equations to accurately portray the dynamics of the full system diminishes and eventually breaks down with increasing interaction strength.

In recent years various strategies have been proposed to overcome the limitations of first-order phase reduction. The most widely employed approach is to include additional variables in directions transverse to the limit cycle which allows for a phase-amplitude description of transient trajectories away from the limit cycle. While some authors have used moving orthonormal coordinate systems to define the amplitudes \cite{Wedgwood2013, Letson2018}, there is a growing body of work that utilises the notion of isostable coordinates \cite{Guillamon2009, Mauroy2013, Wilson2016}. These coordinates represent level sets of the slowest decaying eigenmodes of the Koopman operator \cite{Budisic2012, Mauroy2013} and in the absence of perturbations, have exponential decay to zero within the basin of attraction of the limit cycle at rates given by the Floquet exponents of the linearisation of the flow about the limit cycle. First defined for systems with a stable fixed point \cite{Mauroy2013}, isostable coordinates are now a well established concept for the analysis of perturbed trajectories near periodic orbits \cite{Castejon2013, Perez-Cervera2020, Mauroy2018, Shirasaka2017, Wilson2018} and also for the control of oscillators (see \cite{Wilson2022} for a recent review).

A dimension reduction can be achieved for high-dimensional systems by assuming that isostable coordinates which rapidly decay (associated with negative Floquet exponents of large magnitude) are always zero \cite{Wilson2018, Wilson2019b}. This results in a phase-amplitude reduced description of single oscillator dynamics in response to perturbations which can be taken to arbitrary orders of accuracy in the isostable coordinates \cite{Wilson2019b, Wilson2019, Wilson2020b}. In turn, through high order asymptotic expansion of interactions in terms of the phase and isostable coordinates \cite{Park2021}, it is possible to build network equations describing the evolution of phase and isostable coordinates of each node in a coupled oscillator network. To date, analysis of these network equations has been relatively limited and mainly restricted to networks with just two oscillators \cite{Wilson2019c, Park2021, Ermentrout2019, Coombes2023b}.

Two different but related frameworks for the use of the phase-isostable network differential equations for the study of oscillator network dynamics have emerged. The most frequently employed approach has been to use the isostable dynamics to derive a phase equation to higher order in the coupling strength. The isostable dynamics are assumed to be slaved to the phases so that at each order in the coupling strength the isostables are described by a linear ordinary differential equation \cite{Wilson2019c, Park2021, Park2023, Mau2023}. This gives a higher-order description of the phase dynamics when the solution is substituted into the equation for the time evolution of the phase. The higher-order phase reduced equations for networks of more than two nodes contain non-pairwise phase interactions (i.e., terms involving the phases of three or more oscillators, noted in \cite{Wilson2019c, Park2023, Mau2023}, but overlooked in \cite{Park2021}), despite the interactions between the nonlinear oscillators being pairwise. This is a typical result for higher-order phase reduced network equations \cite{Battison2020}. Similar non-pairwise phase interactions were identified in \cite{Leon2019} where an alternative strategy was used to compute phase reduced equations to second and third order in the coupling strength for the mean-field complex Ginzburg-Landau equation relying on the explicit expression for the isochrons (which can be obtained in this case). 

The alternative framework retains both phase and one isostable coordinate for each node in the coupled oscillator network. One work to have adopted this approach is that of Ermentrout and Wilson \cite{Ermentrout2019}. Similarly to \cite{Wilson2019c, Park2021}, they consider as an example a pair of synaptically coupled thalamic neurons. All three works \cite{Wilson2019c, Park2021, Ermentrout2019} show that their approaches can reveal bifurcations leading to stable phase-locked states at higher coupling strength as observed in full system simulations. Coombes et al., \cite{Coombes2023b} also use network equations for the dynamics of phases and isostable coordinates to consider two linearly coupled planar piecewise-linear caricatures of the Morris-Lecar neuron model. Accounting for the non-smoothness of the dynamics using results from \cite{Wilson2019}, a bifurcation diagram for the dynamics in the phase-isostable framework under variation of coupling strength reveals restabilisation of synchrony, the existence of stable phase-locked states other than synchrony and antisynchrony and also stable quasiperiodic states, none of which can be captured using first-order phase reduction of the dynamics (which predicts only stable antisynchrony and unstable synchrony for weak coupling). All of the dynamics revealed by the phase-isostable network analysis are shown to be in qualitative agreement with numerical simulations of the full system, however comparison with results using the more accurate master stability function \cite{Pecora1998} shows that the phase-isostable approach underestimates the value of the coupling strength at which synchrony restabilises \cite{Coombes2023b}.

Previous work therefore indicates that for networks of two oscillators both phase-amplitude network equations and higher-order phase reduction based on phase-amplitude network interactions have the capability to reveal dynamics which cannot be explained using standard phase reduction. It remains an open problem to fully explore the capabilities of both approaches and draw comparisons between the two in terms of the accuracy with which they can capture the existence and stability of network states and to extend analysis of the approaches to larger networks. 

In this paper we go some way to address these challenges. Following discussion of the necessary background on isostable coordinates in section \ref{sec:background}, we derive the phase-isostable network differential equations and also the corresponding higher-order phase-reduced equations for a network of $N$ identical oscillators in sections \ref{sec:background} and \ref{sec:compare} respectively. Using first-order averaging we obtain a system of phase-isostable network differential equations linear in the isostable coordinates and with six pairwise phase difference interaction functions. These equations are used to determine conditions for the existence and stability of various phase-locked states in section \ref{sec:phaselocked}. We are then able in section \ref{sec:compare} to compare results regarding the location of stability boundaries for full synchrony and the asynchronous splay state arising from each reduction method for the mean--field complex Ginzburg--Landau equation where it is also possible to derive the exact location of the stability boundaries for the full system. We observe that retaining the isostable coordinate through using the phase-isostable network equations gives the greatest accuracy in approximating the qualitative bifurcation structure. 

We therefore proceed to use the phase-isostable network equations to approximate the dynamics of both small and large networks of neural oscillators in section \ref{sec:ML}. For a two node network of linearly coupled Morris-Lecar neurons \cite{Morris1981} we compare the bifurcation diagrams for the full model with that for the phase-isostable network equations. We observe that the phase-isostable description qualitatively agrees with the full model in capturing bifurcations of branches of synchronous and antisynchronous solutions in addition to the existence of stable phase-locked states. This marks a vast improvement on the descriptive power of the first-order phase reduction which cannot capture any of these phenomena. We also investigate a larger network of many globally linearly coupled Morris-Lecar neurons where we find that as for the full model considered in \cite{Han1995} there is an interval in the coupling strength where synchrony and the uniform incoherent (splay) state are both unstable. In this interval we observe stable cluster states in both the full dynamics and the phase-amplitude reduced equations. Finally in section \ref{sec:conclusions} we discuss the limitations of the framework as well as possible extensions of our work to further refine the accuracy of the description of the network dynamics. We also highlight other phenomena which may be revealed in driven phase oscillator dynamics through using phase-amplitude coordinates.

\section{Phase and isostable coordinates and reductions}\label{sec:background}

We consider the dynamics of $N$ identical pairwise coupled oscillators
\begin{align}\label{eq:fullnetwork}
\FD{\vec{x}_i}{t} = \vec{F}(\vec{x}_i) + \epsilon \sum_{j=1}^N w_{ij} \vec{G}(\vec{x}_i, \vec{x}_j), \quad i=1, \ldots N,
\end{align}
where $\vec{x}_i \in \mathbb{R}^n$. The coupling function $\vec{G}: \mathbb{R}^n \times \mathbb{R}^n\to \mathbb{R}^n$ describes the pairwise interactions which take the same form for each pair of oscillators. The overall strength of interactions is given by $\epsilon$ while $w_{ij}$ modulates the weight of connectivity from node $j$ to node $i$, effectively describing the topology of the oscillator network. The uncoupled dynamics (for $\epsilon=0$) are described by $\vec{F}: \mathbb{R}^n \to \mathbb{R}^n$, where we assume that
\begin{align}\label{eq:uncoupled}
\FD{\vec{x}}{t}= \vec{F}(\vec{x})
\end{align}
has a T-periodic hyperbolic limit cycle, $\gamma$ and we denote points on $\gamma$ by $\vec{x}^\gamma(t)$. Let $\phi:\mathbb{R} \times \mathbb{R}^n \to \mathbb{R}^n$ denote the flow induced by \eqref{eq:uncoupled} so that $\phi(t, \vec{x}_0)$ represents the solution of \eqref{eq:uncoupled} with initial condition $\vec{x}(0)= \vec{x}_0$. Then the monodromy matrix $M(\vec{x}^\gamma(t_*)) = \partial \phi(T, \vec{x})/ \partial \vec{x} |_{\vec{x}=\vec{x}^\gamma(t_*)}$ (time--$T$ flow linearised about a point $\vec{x}^\gamma(t_*)$ on the orbit) has eigenvalues $\lambda_m$, $m=0, \ldots, n-1$ called the Floquet (or characteristic) multipliers and we assume that $1=\lambda_0 > |\lambda_1| \geq \ldots \geq |\lambda_{n-1}|$ so that the limit cycle $\gamma$ is stable. For simplicity we further assume that the Floquet multipliers $\lambda_m$ are positive, real and simple (and see \cite{Wilson2019} for discussion on how to proceed when this is not the case). The non-zero Floquet exponents are then $\kappa_m= \log(\lambda_m)/T$, $m=1, \ldots, n-1$ which are all real and negative.

\subsection{First-order phase reduction and weakly coupled phase oscillators}

The periodic orbit $\gamma$ of \eqref{eq:uncoupled} can be parameterised by a phase $\theta \in [0, 2\pi)$ with an arbitrary point $\vec{x}^\gamma_0 \in \gamma$ having $\theta(\vec{x}^\gamma_0)=0$ and $\theta(\phi(t, \vec{x}^\gamma_0))= \omega t$ where $\omega=2\pi/T$. The notion of phase can be extended to points $\vec{x}_*$ in the basin of attraction of the limit cycle $\mathcal{B}(\gamma)$ by defining the asymptotic phase as the unique $\theta(\vec{x}_*) \in [0, 2\pi)$ satisfying
 \begin{align}
\lim_{t\to\infty} \left|\phi(t, \vec{x}_*) - \phi(t + \theta(\vec{x}_*)/\omega, \vec{x}^\gamma_0)\right|=0.
\end{align}
The level sets of $\theta(\vec{x}_*)$ are called isochrons and contain all points in $\mathcal{B}(\gamma)$ with the same asymptotic phase \cite{Winfree1967, Guckenheimer1975, Kuramoto1984, Winfree2001}. The phase of each uncoupled node trajectory then evolves as $\rm{d} \theta/\rm{d} t = \omega$ both on and off the limit cycle. Note that we can then use $\vec{x}^\gamma(t)$ and $\vec{x}^\gamma(\theta)$ to denote points on cycle where time is parameterised as $t=\theta/\omega$. In terms of phase variables, \eqref{eq:fullnetwork} becomes
\begin{align}
\FD{\theta_i}{t} = \left.\PD{\theta_i}{\vec{x}_i}\right|_{\vec{x}_i}\cdot \FD{\vec{x}_i}{t} = \omega + \epsilon \sum_{j=1}^N w_{ij} \left.\PD{\theta_i}{\vec{x}_i}\right|_{\vec{x}_i}\cdot  \vec{G}(\vec{x}_i, \vec{x}_j),
\end{align} where $\left.\PD{\theta_i}{\vec{x}_i}\right|_{\vec{x}_i}$ is the phase response curve which quantifies the effects of a perturbation on the phase of the oscillator. If the coupling is weak, then the dynamics stay in a neighbourhood of the limit cycle and evaluating on the limit cycle we have to first order in $\epsilon$,
\begin{align}\label{eq:phasereduction}
\FD{\theta_i}{t} =  \omega + \epsilon \sum_{j=1}^N w_{ij} Z(\theta_i) \cdot  \vec{G}(\vec{x}^\gamma(\theta_i), \vec{x}^\gamma(\theta_j)),
\end{align} where $Z(\theta) = \left.\PD{\theta}{\vec{x}}\right|_{\vec{x}^\gamma(\theta)} = \nabla_{\vec{x}^\gamma(\theta)}\theta$ is the gradient of the phase variable $\theta$ evaluated on the limit cycle and is known as the infinitesimal phase response curve (iPRC). While $Z$ may be found directly \cite{Winfree2001, Brown2004, Monga2019}, it most commonly computed using the adjoint method \cite{Brown2004, Ermentrout1991, Hoppensteadt1997} whereby $Z$ is the $T$-periodic solution of the adjoint equation
\begin{align}\label{eq:PRCadjoint}
\FD{Z(\theta)}{t} = -J^T Z(\theta),
\end{align} satisfying the normalisation condition $Z(0) \cdot \vec{F}(\vec{x}^\gamma(0))= \omega$.
Here ${}^T$ denotes the transpose and $J:={\rm D}\vec{F}(\vec{x}^\gamma(\theta))$ is the Jacobian of the vector field $\vec{F}$ evaluated on cycle.

A further simplification of \eqref{eq:phasereduction} can be made by transforming the system to a rotating frame through $\theta_i = \phi_i + \omega t$ where, assuming weak coupling, $\phi_i$ will slowly drift. First order averaging \cite{Guckenheimer1990, Ermentrout1991, Hoppensteadt1997} then gives the phase dynamics back in the original variables as the simpler to analyse phase difference equations
\begin{align}\label{eq:phaseaveraged}
\FD{\theta_i}{t} = \omega + \epsilon \sum_{j=1}^{N} w_{ij} H_1(\theta_j-\theta_i), 
\end{align} where \begin{align}
    H_1(\chi) = \frac{1}{2\pi} \int_0^{2\pi} Z(u) \cdot \vec{G}(\vec{x}^\gamma(u), \vec{x}^\gamma(u+\chi))\, {\rm d}u.
\end{align} is known as the phase interaction function. A solution of the averaged equations \eqref{eq:phaseaveraged} is $\epsilon$-close to a solution of the unaveraged equations \eqref{eq:phasereduction} for times of $O(\epsilon^{-1})$. The existence and stability of phase-locked states in the first-order phase regime can be studied using \eqref{eq:phaseaveraged} and specified in terms of the properties of $H_1$ and the choice of network structure given by $W = (w_{ij})$ \cite{Ashwin1992, Kim2000, Watanabe1997, Ermentrout1992, Coombes2023}. In section \ref{sec:backgroundiso} we will show how the theory of isostables can be used to derive averaged phase-isostable network equations  with six interaction functions $H_k$, $k=1, \ldots, 6$ and in section \ref{sec:phaselocked} we investigate the existence and stability of phase locked states for larger magnitude interactions in terms of properties of the $H_k$, and network structure.

Equation \eqref{eq:phasereduction} gives what we will refer to as the first-order phase reduction of the network dynamics. For further details on the phase reduction of oscillator networks see the recent reviews \cite{Kuramoto2019, Pietras2019}. Phase-reduced equations have provided the standard framework for understanding dynamics of weakly coupled oscillator networks for the last four decades and can be extremely instructive. For instance, the Kuramoto model \cite{Kuramoto1984} has simply $H_1(\chi)= \sin(\chi)$ yet is able to capture the basic mechanisms underlying synchronisation in many biological \cite{Cumin2007}, chemical \cite{Forrester2015} and physical \cite{Wiesenfeld1998} oscillator networks. However, the phase-reduction \eqref{eq:phasereduction} is only valid provided the state of each oscillator remains close to its underlying limit cycle. This requires that interactions are weak ($|\epsilon| \ll 1$) and therefore \eqref{eq:phasereduction} cannot be used to describe network dynamics with interactions of larger magnitude which could potentially lead to oscillators spending significant time away from the limit cycle. In order to capture transient dynamics away from the limit cycle we introduce amplitude variables in directions transverse to the limit cycle.

\subsection{Isostables and phase-isostable reduction of network equations}\label{sec:backgroundiso}

In this paper we use the concept of isostables to define coordinates off limit cycle \cite{Guillamon2009, Shirasaka2017, Wilson2018, Wilson2019, Monga2019, Mauroy2013, Wilson2016}. Isostables for the uncoupled system \eqref{eq:uncoupled} identify initial conditions with the same relaxation rate to the limit cycle and therefore approach the limit cycle together \cite{Shirasaka2017}. For each Floquet exponent $\kappa_m$, $m=1, \ldots, n-1$, a set of isostables representing an amplitude degree of freedom can be defined and therefore for any point $\vec{x} \in \mathcal{B}(\gamma)$ we can associate isostable coordinates $\psi_m$, $m=1, \ldots, n-1$. The isostables can be defined as level sets of certain eigenfunctions of the Koopman operator \cite{Mauroy2013, Kvalheim2021, Shirasaka2017, Mauroy2018} and a constructive definition may be given for the slowest decaying isostable $\psi_1$ \cite{Wilson2018, Wilson2020b}: Denote by $w^T$ the left eigenvector of $M(\vec{x}_0^\gamma)$ associated with the eigenvalue $\lambda_1 = \exp(\kappa_1 T)$. The corresponding right eigenvector is $v$ with $\|v\| = \sqrt{v \cdot v} =1$ and $w^T v= 1$. Then
\begin{align}\label{eq:isostablesWilson}
 \psi_1(\vec{x}) = \lim_{k\to \infty}\left( w^T(\phi(t_{\theta=0}^k , \vec{x}) - \vec{x}_0^\gamma) \exp(-\kappa_1 t_{\theta=0}^k) \right)
\end{align}
where $t_{\theta=0}^k$ is the time of the $k$th transversal of the $\theta=0$ isochron. The isostable coordinates are amplitudes which satisfy $\dot{\psi}_m = \kappa_m \psi_m$ in the absence of coupling or perturbations.

An alternative perspective is taken by the authors of \cite{Guillamon2009, Castejon2013, Perez-Cervera2020} who note that for each point $\vec{x} \in \mathcal{B}(\gamma)$ there is a coordinate transform given by an analytic map $K$ such that $\vec{x} = K(\theta, \psi_1, \ldots \psi_{n-1})$, and determine the map $K$ using a parameterization method.

For oscillator $i$ in the coupled network \eqref{eq:fullnetwork} we can now describe the phase and isostable dynamics by
\begin{subequations}\label{eqs:pinetwork}
\begin{align}\label{eq:phaseiso}
 \FD{\theta_i}{t} & = \omega + \epsilon \sum_{j=1}^N w_{ij} \left.\PD{\theta_i}{\vec{x}_i}\right|_{\vec{x}_i} \cdot  \vec{G}(\vec{x}_i, \vec{x}_j),  \\ \label{eq:phaseiso2}
 \FD{(\psi_i)_m}{t} & = \kappa_m (\psi_i)_m + \epsilon \sum_{j=1}^N w_{ij} \left.\PD{(\psi_i)_m}{\vec{x}_i}\right|_{\vec{x}_i}\cdot  \vec{G}(\vec{x}_i, \vec{x}_j),\\ & \qquad \qquad m=1, \ldots, n-1 \notag,
\end{align}
\end{subequations}
where the gradients $\left.\PD{\theta_i}{\vec{x}_i}\right|_{\vec{x}_i}$ and  $\left.\PD{(\psi_i)_m}{\vec{x}_i}\right|_{\vec{x}_i}$ are the phase and isostable response curves respectively that quantify the  effects of a perturbation on the phase and amplitude coordinates of the oscillator. We now make the simplifying assumption that $\kappa_2, \ldots, \kappa_{n-1}$ are all sufficiently close to zero that perturbations in the directions of isostables $(\psi_i)_2, \ldots, (\psi_i)_{n-1}$ may be ignored. That is we consider a single isostable coordinate $(\psi_i)_1:= \psi_i$ for each oscillator with corresponding Floquet exponent $\kappa_1 := \kappa$. This assumption is similar to that made in previous work \cite{Wilson2018,Wilson2019c,Park2021}. We have therefore (for node dynamics of dimension $n>2$) made a reduction in order from $Nn$ equations to $2N$ equations for the network dynamics, retaining only the most slowly decaying isostable dynamics for each oscillator.

Following \cite{Wilson2020b, Park2021} we can take asymptotic expansions of solutions away from the limit cycle and also the phase and isostable response curves about the limit cycle to arbitrary order in the $O(\epsilon)$ isostable coordinate. That is we can write
\begin{subequations}\label{eqs:expansions}
\begin{align}\label{eq:expg}
\vec{x}(\theta, \psi) & = \vec{x}^\gamma(\theta) + \Delta \vec{x}(\theta, \psi) = \vec{x}^\gamma(\theta) + \sum_{k=1}^\infty \psi^kg^{(k)}(\theta),\\ \label{eq:expZ}
\left.\PD{\theta}{\vec{x}}\right|_{\vec{x}} &= \mathcal{Z}(\theta, \psi) = Z^{(0)}(\theta) + \sum_{k=1}^\infty \psi^k Z^{(k)}(\theta), \\ \label{eq:expI}
\left.\PD{\psi}{\vec{x}}\right|_{\vec{x}} &= \mathcal{I}(\theta, \psi)= I^{(0)}(\theta) + \sum_{k=1}^\infty \psi^k I^{(k)}(\theta),
\end{align}
\end{subequations} where $g^{(1)}(\theta)$ is the Floquet eigenfunction (right eigenvector of $M(\vec{x}^\gamma(\theta))$ associated with eigenvalue $e^{\kappa T}$) and $g^{(k)}(\theta)$, $k>1$ are higher order analogues. The gradient of the phase and amplitude coordinates evaluated on the limit cycle at phase $\theta$ are the iPRC $Z^{(0)}(\theta)$ and the infinitesimal isostable response curve (iIRC) $I^{(0)}(\theta)$ respectively. The terms $Z^{(k)}(\theta)$ and $I^{(k)}(\theta)$, $k>0$, are higher order correction terms to $Z^{(0)}(\theta)$ and $I^{(0)}(\theta)$. Note that this generalises the linear order expansions developed in \cite{Wilson2018,Wilson2019b,Wilson2019}. All of the $T$-periodic vector functions $g^{(k)}$, $Z^{(k)}$, $I^{(k)}$ can be computed using appropriately normalised adjoint equations \cite{Wilson2020b}:
\begin{subequations}\label{eqs:adjoints}\begin{align}
\FD{g^{(k)}}{t} &= (J - k\kappa I_n)g^{(k)} + \alpha^{(k)}
 \label{eq:eigODE} \\
\FD{{Z}^{(k)}}{t} &= -(J^{T}+k\kappa I_n){Z}^{(k)} - \sum_{i=1}^{n} \sum_{j=0}^{k-1} {\rm e}_{i}^{T}{Z}^{(j)}b_{i}^{(k-j)}, \label{eq:PRCODE} \\
\FD{{I}^{(k)}}{t} &= -(J^{T}+(k-1)\kappa I_n){I}^{(k)} - \sum_{i=1}^{n} \sum_{j=0}^{k-1} {\rm e}_{i}^{T}{I}^{(j)}b_{i}^{(k-j)}, \label{eq:IRCODE}
\end{align}\end{subequations}
where $I_n$ is the $n\times n$ identity matrix, ${\rm e}_{i}$ are the unit basis vectors and $\alpha^{(k)}= [ \alpha_1^{(k)}\  \ldots\  \alpha_n^{(k)}]^T$ and $b_{i}^{(k-j)}$ are the vectors defined in \eqref{eq:alpha} and \eqref{eq:bvectors} respectively.  See \cite{Wilson2020b} for the derivation of these equations for all $N$ isostable coordinates and Appendix \ref{sec:AppendixA} for the derivation of the compact notation \eqref{eqs:adjoints} when a single isostable coordinate is considered. The appropriate normalisations for \eqref{eqs:adjoints} are also given in Appendix \ref{sec:AppendixA}. In general the solutions of the hierarchy of equations \eqref{eqs:adjoints} must be computed numerically \cite{Monga2019, Park2021}.

In order to obtain the network equations \eqref{eqs:pinetwork} in terms of phase and isostable coordinates, we also require the appropriate expansion for the coupling function $\vec{G}: \mathbb{R}^n \times \mathbb{R}^n\to \mathbb{R}^n$ which can be expressed to arbitrary order following \cite{Park2021}. Here we will require only the expansion up to and including second order derivatives of $\vec{G}$ and therefore we may use more convenient notation than that in \cite{Park2021}. For $\vec{X} = [\vec{x}_i^T,\vec{x}_j^T]^T \in \mathbb{R}^{2n}$ we write $\vec{G}(\vec{X}) = [ G_1(\vec{X}) \ldots G_n(\vec{X})]^T$ where $G_q(\vec{X}) \in \mathbb{R}$. Defining $\vec{X}^\gamma = [\vec{x}^\gamma(\theta_i)^T,\vec{x}^\gamma(\theta_j)^T]^T$ and $\Delta \vec{X} = [\Delta \vec{x}_i^T,\Delta \vec{x}_j^T]^T$, the Taylor expansion of $\vec{G}$ can be expressed to second order as
\begin{align}\label{eq:expG}
\vec{G}(\vec{X}^\gamma + \Delta \vec{X}) = \vec{G}(\vec{X}^\gamma) +  \vec{DG}(\vec{X}^\gamma)\Delta \vec{X} + \frac{1}{2} \begin{bmatrix}
\Delta \vec{X}^T \vec{H}_1 \Delta \vec{X} \\
\vdots \\
\Delta \vec{X}^T \vec{H}_n \Delta \vec{X} \\
\end{bmatrix} + \cdots
\end{align}
where $\vec{D} \vec{G}$ is the Jacobian of $\vec{G}$ with respect to $\vec{X}$ and $\vec{H}_q$ is the Hessian matrix of second order derivatives of $G_q$ where and all derivatives are evaluated at $\vec{X}^\gamma$. Since
\begin{align}
  \vec{DG}(\vec{X}^\gamma) = \left[ J_1\ J_2 \right], \qquad \vec{H}_q = \begin{bmatrix}
 \vec{H}_q^{11} & \vec{H}_q^{12}\\ \vec{H}_q^{21} & \vec{H}_q^{22}
 \end{bmatrix}
\end{align} where $J_k= \PD{\vec{G}}{\vec{x}_k}$ is the Jacobian of $\vec{G}(\vec{x}_i, \vec{x}_j)$ with respect to its $k$th argument and $\vec{H}_q^{kl} = \PD{}{\vec{x}_l^T}\left( \PD{G_q}{\vec{x}_k}\right)$ we can rewrite \eqref{eq:expG} as
\begin{align}
\begin{split}
\vec{G}(\vec{x}_i, \vec{x}_j) &= \vec{G}(\vec{x}^\gamma(\theta_i),\vec{x}^\gamma(\theta_j)) + J_1 \Delta \vec{x}_i + J_2 \Delta \vec{x}_j \\  &\ \ + \frac{1}{2} \begin{bmatrix}\Delta \vec{x}_i^T \vec{H}_1^{11} \Delta \vec{x}_i  + \Delta \vec{x}_j^T \vec{H}_1^{22} \Delta \vec{x}_j \\ \vdots\\ \Delta \vec{x}_i^T \vec{H}_n^{11} \Delta \vec{x}_i  + \Delta \vec{x}_j^T \vec{H}_n^{22} \Delta \vec{x}_j \end{bmatrix} \\ &\ \ + \frac{1}{2} \begin{bmatrix} \Delta \vec{x}_i^T \vec{H}_1^{12} \Delta \vec{x}_j + \Delta \vec{x}_j^T \vec{H}_1^{21} \Delta \vec{x}_i  \\ \vdots\\  \Delta \vec{x}_i^T \vec{H}_n^{12} \Delta \vec{x}_j + \Delta \vec{x}_j^T \vec{H}_n^{21} \Delta \vec{x}_i \end{bmatrix} + \cdots.
\end{split}
\end{align} Using the expansions of $\Delta \vec{x}_i$ and $\Delta \vec{x}_j$ from \eqref{eq:expg} and collecting terms up to quadratic order in $\psi$ we have
\begin{align}\label{eq:couplingexp} \begin{split}
    \vec{G}(\theta_i, \psi_i, \theta_j, \psi_j) =& \vec{G}(\vec{x}^\gamma(\theta_i),\vec{x}^\gamma(\theta_j)) + \psi_i J_1 g^{(1)}(\theta_i) \\  & +  \psi_j J_2 g^{(1)}(\theta_j)  + \psi_i^2 K_1(\theta_i, \theta_j) \\  & +  \psi_j^2 K_2(\theta_i, \theta_j)  + \psi_i \psi_j L(\theta_i, \theta_j) + \cdots ,\end{split}
\end{align} where
\begin{subequations}
\begin{align}
K_1(\theta_i, \theta_j) &= J_1 g^{(2)}(\theta_i) + \frac{1}{2} \begin{bmatrix} g^{(1)}(\theta_i)^T \vec{H}^{11}_1 g^{(1)}(\theta_i)\\ \vdots \\  g^{(1)}(\theta_i)^T \vec{H}^{11}_n g^{(1)}(\theta_i)\end{bmatrix}\\
K_2(\theta_i, \theta_j) &= J_2 g^{(2)}(\theta_j) + \frac{1}{2} \begin{bmatrix} g^{(1)}(\theta_j)^T \vec{H}^{22}_1 g^{(1)}(\theta_j)\\ \vdots \\  g^{(1)}(\theta_j)^T \vec{H}^{22}_n g^{(1)}(\theta_j)\end{bmatrix}\\
L(\theta_i, \theta_j) &= \frac{1}{2} \begin{bmatrix} g^{(1)}(\theta_i)^T \vec{H}^{12}_1 g^{(1)}(\theta_j) + g^{(1)}(\theta_j)^T \vec{H}^{21}_1 g^{(1)}(\theta_i)\\ \vdots \\  g^{(1)}(\theta_i)^T \vec{H}^{12}_n g^{(1)}(\theta_j)+ g^{(1)}(\theta_j)^T \vec{H}^{21}_n g^{(1)}(\theta_i)\end{bmatrix}.
\end{align}
\end{subequations}

We now have all of the expansions required to express \eqref{eq:fullnetwork} in terms of phase and isostable coordinates. Including terms up to $O(\epsilon^2)$ the phase-amplitude network equations are
\begin{subequations}\label{eqs:pinetwork1}
    \begin{align}\label{eq:phase}
    \begin{split}
 \FD{\theta_i}{t} & = \omega + \epsilon \sum_{j=1}^N w_{ij} \Bigl[h_1(\theta_i, \theta_j) + \psi_i h_2(\theta_i, \theta_j) \\  & \qquad \qquad \qquad \qquad \qquad \quad+ \psi_j h_3(\theta_i, \theta_j)\Bigr], \end{split} \\ \label{eq:iso}
 \begin{split}
 \FD{\psi_i}{t} & = \kappa \psi_i + \epsilon \sum_{j=1}^N w_{ij} \Bigl[h_4(\theta_i, \theta_j) + \psi_i h_5(\theta_i, \theta_j) \\  & \qquad \qquad \qquad \qquad \qquad \quad  + \psi_j h_6(\theta_i, \theta_j)\Bigr],\end{split}
\end{align}
\end{subequations} where the six interaction functions $h_k(\theta_j, \theta_j)$ are given by
\begin{subequations}\label{eqs:hs}
\begin{align}\label{eq:h1}
h_1(\theta_i, \theta_j) &= Z^{(0)}(\theta_i)\cdot \vec{G}(\vec{x}^\gamma(\theta_i),\vec{x}^\gamma(\theta_j)), \\
h_2(\theta_i, \theta_j) &=  Z^{(0)}(\theta_i) \cdot J_1 g^{(1)}(\theta_i) + Z^{(1)}(\theta_i)\cdot \vec{G}(\vec{x}^\gamma(\theta_i),\vec{x}^\gamma(\theta_j)), \\
h_3(\theta_i, \theta_j) &= Z^{(0)}(\theta_i)\cdot J_2 g^{(1)}(\theta_j),\\ \label{eq:h4}
h_4(\theta_i, \theta_j) &= I^{(0)}(\theta_i)\cdot \vec{G}(\vec{x}^\gamma(\theta_i),\vec{x}^\gamma(\theta_j)),\\\label{eq:h5}
h_5(\theta_i, \theta_j) &=  I^{(0)}(\theta_i) \cdot J_1 g^{(1)}(\theta_i) + I^{(1)}(\theta_i)\cdot \vec{G}(\vec{x}^\gamma(\theta_i),\vec{x}^\gamma(\theta_j)), \\ \label{eq:h6}
h_6(\theta_i, \theta_j) &= I^{(0)}(\theta_i)\cdot J_2 g^{(1)}(\theta_j).
\end{align}\end{subequations}
The expansion of the network equations may be computed to higher order, but note that to $O(\epsilon^m)$ there are $m(m+1)$ interaction functions. Here we explore the capabilities of the $O(\epsilon^2)$ expansion \eqref{eqs:pinetwork1} to describe network behaviours and discuss the limitations imposed by the truncation at this order in section \ref{sec:conclusions}. The network equations (and therefore analysis) may be simplified to a phase difference system by first-order averaging \cite{Guckenheimer1990, Sanders2007, Ermentrout2019} yielding
\begin{subequations}\label{eqs:piaveraged}
    \begin{align} \begin{split}
        \FD{\theta_i}{t} =& \omega + \epsilon\sum_{j=1}^{N} w_{ij}\Bigl[H_{1}(\theta_j-\theta_i)+\psi_iH_2(\theta_j-\theta_i) \\ & \qquad \qquad \qquad \qquad \qquad  +\psi_jH_3(\theta_j-\theta_i)\Bigr], \label{eq:avphase}
    \end{split} \\ \begin{split}
\FD{\psi_i}{t} =& \kappa\psi_i + \epsilon\sum_{j=1}^{N} w_{ij}\Bigl[H_{4}(\theta_j-\theta_i)+\psi_i H_5(\theta_j-\theta_i) \\ & \qquad \qquad \qquad \qquad \qquad  +\psi_j H_6(\theta_j-\theta_i)\Bigr], \label{eq:aviso} \end{split}
\end{align}
\end{subequations}
where $H_k$ are the $2\pi$-periodic functions
\begin{align}
H_k(\chi) = \frac{1}{2\pi} \int_{0}^{2\pi} h_k(u,u+\chi) \,{\rm d}u. \label{eq:Hk}
\end{align} In the weak coupling limit where the isostable coordinates may be assumed to be zero, we recover the first-order phase reduced equation \eqref{eq:phaseaveraged}. The equations \eqref{eqs:piaveraged} generalise to networks of arbitrary size and structure those given in \cite{Ermentrout2019} and contain higher order terms than those derived in \cite{Wilson2019c}. Analysis of existence and stability of phase-locked states in \eqref{eqs:piaveraged} is investigated next in Section \ref{sec:phaselocked}. Note that since \eqref{eqs:piaveraged} are averaged equations, the solutions of \eqref{eqs:pinetwork1} and \eqref{eqs:piaveraged} are $\epsilon$-close for times of $O(\epsilon^{-1})$.

\section{Phase-locked states in phase-amplitude network equations} \label{sec:phaselocked}

Having defined a system of phase-amplitude network equations \eqref{eqs:piaveraged}, we now investigate conditions for the existence and stability of certain phase-locked states within the averaged equations, generalising well known results for phase-reduced networks \cite{Ashwin1992, Kim2000, Watanabe1997, Ashwin2016, Coombes2023}.

A $1:1$ phase-locked state in the network of $N$ identical oscillators \eqref{eqs:piaveraged} is defined by $\theta_i = \phi_i + \Omega t$ where the $\phi_i$ are constant phase lags and $\Omega$ is the collective frequency of the oscillators. We denote such a state by $\Phi = (\phi_1, \ldots, \phi_N)$. In phase-reduced systems all node orbits are assumed to coincide with the stable limit cycle, $\gamma$ of the node dynamics. By including information about the dynamics off-cycle through the isostable coordinates $\psi_i$ we allow each node to have a different trajectory in the phase-isostable phase space. That is, for a solution $(\theta_1, \ldots, \theta_N, \psi_1, \ldots, \psi_N)$ of \eqref{eqs:piaveraged}, the projection $(\theta_i, \psi_i)= (\phi_i + \Omega t, \psi_i(t))$ can be different for each node. Substituting into the $O(\epsilon^2)$ truncated averaged equations \eqref{eqs:piaveraged} we have for $i = 1, \ldots, N$,
\begin{subequations}\label{eqs:phaselocked1}
\begin{align}
\begin{split}\label{eq:pl1}
&\Omega = \omega + \epsilon \sum_{j=1}^N w_{ij} \Bigl[ H_1(\phi_j - \phi_i) + \psi_i(t) H_2(\phi_j-\phi_i)  \\ & \qquad\qquad \qquad \qquad \qquad + \psi_j(t) H_3(\phi_j-\phi_i)\Bigr],\end{split}\\ \label{eq:pl2}\begin{split}
&\FD{\psi_i}{t} = \kappa \psi_i(t) + \epsilon \sum_{j=1}^N w_{ij} \Bigl[ H_4(\phi_j - \phi_i) \\ &  \qquad \qquad +\psi_i(t) H_5(\phi_j-\phi_i)  + \psi_j(t) H_6(\phi_j-\phi_i)\Bigr].\end{split}
\end{align} 
\end{subequations} For a given fixed set of relative phases $\Phi$ and connectivity $W=(w_{ij})$, equations \eqref{eq:pl1} are $N$ linear equations from which we determine that $\psi_i(t)$ are constants depending on $\Omega$. That is, any phase-locked solution of \eqref{eqs:piaveraged} has node orbits which coincide with an isostable. Note that this result is particular to the truncation of the network equations \eqref{eqs:piaveraged} at $O(\epsilon^2)$ which includes only linear terms in the isostable coordinates and would not be expected to hold for higher order truncation.

Denoting the constant isostable coordinates $\psi_i(t) = {\Psi}_i$ we find from \eqref{eqs:phaselocked1} that a phase locked state with relative phases $\Phi$ exists with collective frequency $\Omega$ only if 
\begin{align}(\Omega- \omega) \vec{1}_N = \epsilon(\vec{p} - P Q^{-1} \vec{q})\label{eq:existence}\end{align}
where $\vec{1}_N$ is the column vector of 1s, $\vec{p} = (p_1, \ldots, p_N)^T$, $\vec{q} = (q_1, \ldots, q_N)^T$  with 
\[p_i = \sum_{j=1}^N w_{ij} H_1(\phi_j-\phi_i),\quad  q_i = \sum_{j=1}^N w_{ij} H_4(\phi_j-\phi_i), \] 
and $P$ and $Q$ are matrices with entries
\begin{subequations}
\begin{align}
    P_{ij} &= -w_{ij}H_3(\phi_j-\phi_i) - \delta_{ij} \sum_{k=1}^N w_{ik} H_2(\phi_k-\phi_i) ,\\
    Q_{ij} &= - w_{ij}H_6(\phi_j-\phi_i) - \delta_{ij} \left(\sum_{k=1}^N w_{ik}H_5(\phi_k-\phi_i) +\kappa/\epsilon\right).
\end{align} 
\end{subequations} In this case the constant isostable values for the nodes are given by $\Psi = (\Psi_1, \ldots, \Psi_N)^T = Q^{-1}\vec{q}$. 

We now consider conditions for the stability of a general phase-locked state $\Phi= (\phi_1, \ldots, \phi_N)$ in the $O(\epsilon^2)$ truncated phase-isostable network equations \eqref{eqs:piaveraged}. We linearise about the phase-locked solution by setting
\begin{align}
\theta_i(t) &= \phi_i + \Omega t + \Delta\theta_i(t),\\
\psi_i(t) &= {\Psi}_i + \Delta\psi_i(t),
\end{align} where $\Delta\theta_i(t)$ and $\Delta\psi_i(t)$ are perturbations along and transverse to the periodic orbit $\Gamma$ respectively. The linearisation of \eqref{eqs:piaveraged} gives
\begin{subequations}
    \begin{align}
\begin{split}
 \FD{\Delta\theta_i}{t} \simeq \ & \epsilon\sum_{j=1}^N w_{ij}\Bigl[ H'_{1}(\phi_j - \phi_i)(\Delta\theta_j - \Delta\theta_i)\\
&+   \left({\Psi}_i H'_{2}(\phi_j - \phi_i) +   {\Psi}_j H'_{3}(\phi_j - \phi_i)\right) (\Delta\theta_j - \Delta\theta_i) 
\\
& + \Delta\psi_i H_{2}(\phi_j - \phi_i)
+  \Delta\psi_j H_{3}(\phi_j - \phi_i)\Bigr],
\end{split} \label{eq:ExpandPhase}\\
\begin{split}
\FD{\Delta\psi_i}{t} \simeq \ & \kappa \Delta\psi_i + \epsilon\sum_{j=1}^N w_{ij}\Bigl[ H'_{4}(\phi_j - \phi_i)(\Delta\theta_j - \Delta\theta_i)\\
&+ \left( {\Psi}_i H'_{5}(\phi_j - \phi_i) + {\Psi}_j H'_{6}(\phi_j - \phi_i)\right)(\Delta\theta_j - \Delta\theta_i)\\
&+  \Delta\psi_i H_{5}(\phi_j - \phi_i)
+ \Delta\psi_j H_{6}(\phi_j - \phi_i)\Bigr],
\end{split} \label{eq:ExpandIsostable}
\end{align} \end{subequations}
where $H_k'(\chi) = {\rm d} H_k/ {\rm d} \chi$. Therefore we obtain
\begin{subequations}
\begin{align}
\FD{\Delta\theta_i}{t}&=\sum_{j=1}^N \left[\mathcal{H}_{ij}^{(1)}(\Phi)\Delta\theta_j + \mathcal{H}_{ij}^{(2)}(\Phi)\Delta\psi_j\right], \label{eq:PhaseMatrix}\\
\FD{\Delta\psi_i}{t}&=\sum_{j=1}^N \left[\mathcal{H}_{ij}^{(3)}(\Phi)\Delta\theta_j + \mathcal{H}_{ij}^{(4)}(\Phi)\Delta\psi_j\right], \label{eq:IsostableMatrix}
\end{align}\end{subequations} and hence the Jacobian is the $2N \times 2N$ block matrix
\begin{align}\label{eq:pljac} \mathcal{J} =
\left[\begin{array}{@{}c|c@{}}
  \mathcal{H}^{(1)}(\Phi) & \mathcal{H}^{(2)}(\Phi) \\
\hline
  \mathcal{H}^{(3)}(\Phi) & \mathcal{H}^{(4)}(\Phi)
\end{array}\right]
\end{align}
where the $N \times N$ matrices $\mathcal{H}^{(q)}(\Phi)$ are given in Appendix~\ref{sec:AppendixB}. The phase-locked state $\Phi$ is stable if the eigenvalues of $\mathcal{J}$ all have negative real part with the exception of the eigenvalue which is forced to be zero corresponding to perturbations along the periodic orbit. In general the $\mathcal{H}^{(q)}$ are not in graph-Laplacian form as they are for the phase-only case \cite{Ashwin2016, Coombes2023} and therefore the eigenvalues of $\mathcal{J}$ cannot be related directly to the eigenvalues of $W=(w_{ij})$.  Nonetheless we may still derive conditions required for the existence and stability of particular phase-locked solutions. 

\subsection{Phase-locked states with a common orbit}

We first consider phase-locked states where all oscillators have a common orbit with isostable coordinate $\Psi_i = \Psi$ for all $i=1, \ldots, N$. Such states include full synchrony ($\phi_i= 0$ for all $i$), the splay (or uniform incoherent) state ($\phi_i = 2\pi i/N$) and balanced clusters where the $N= Mq$ nodes lie in $M$ clusters, each containing $q$ nodes and phase difference between the clusters is $2\pi/M$. In Section \ref{sec:twoclusters} we consider unbalanced two-cluster states where each cluster is a different size and has a different orbit.

\subsubsection{Full synchrony}\label{sec:sync}

Synchrony is the most widely studied example of a $1:1$ phase-locked state. We extend established results regarding the stability of synchrony in weakly coupled (phase-reduced) networks by using the phase-isostable network equations \eqref{eqs:piaveraged}. Here all nodes share a common phase $\theta_1(t) = \ldots = \theta_N(t)$ with $\dot{\theta_i} = \Omega$ and orbit with isostable coordinate $\Psi$. Without loss of generality, let $\phi_i = 0$. Then from \eqref{eqs:phaselocked1}
\begin{subequations}\label{eqs:sync}
\begin{align}\label{eq:sync1}
\Omega &= \omega + \epsilon\sum_{j=1}^N w_{ij}\left[H_1(0) + \Psi\left(H_2(0)+H_3(0)\right)\right]\\ \label{eq:sync2}
0 &= \kappa\Psi + \epsilon \sum_{j=1}^N w_{ij}\left[H_4(0) + \Psi\left(H_5(0)+H_6(0)\right)\right],
\end{align}    
\end{subequations} for all $i= 1,\ldots, N$. Denoting the row sum $\sum_{j=1}^N w_{ij}= c_i$, we find that
\begin{align}\label{eq:synciso}
\Psi = & - \frac{\epsilon c_i H_4(0)}{\kappa + \epsilon c_i (H_5(0) + H_6(0))} \\ \label{eq:syncomega}
\Omega = & \omega + \epsilon c_i H_1(0) - \epsilon^2 c_i^2 \frac{ H_4(0)(H_2(0)+ H_3(0))}{\kappa + \epsilon c_i (H_5(0) + H_6(0))}.
\end{align} The values of $\Omega$ and $\Psi$ must be identical for all $i$ to guarantee the existence of the synchronous state. There are two ways this may be achieved.
\begin{enumerate}
\item[a.] If the row sums $c_i$ are independent of $i$ (i.e. $c_i \equiv c$ for all $i$). This is the case for global coupling where $w_{ij}= 1/N$ which we will consider further later.
\item[b.] If $H_1(0)=0$ and $H_4(0)=0$. 
\end{enumerate}
We consider the conditions required for linear stability in both cases. 

\paragraph{Connectivity matrix with constant row sums}
Taking $c_i =c$ for all $i$, \eqref{eq:Hij1}--\eqref{eq:Hij4} give the Jacobian for the synchronous state as
\begin{align}\label{eq:Jconstrows}\begin{split}
    \mathcal{J} =& -\epsilon\begin{bmatrix}
(H_1'(0) + {\Psi}(H_2'(0) + H_3'(0)) &  H_3(0) \\
(H_4'(0) + {\Psi}(H_5'(0) + H_6'(0)) &  H_6(0)
\end{bmatrix} \otimes \mathcal{L} \\ &+ \begin{bmatrix}
0 & \epsilon c(H_2(0)+ H_3(0)) \\
0 & \kappa + \epsilon c (H_5(0) + H_6(0))
\end{bmatrix} \otimes I_N,\end{split}
\end{align} where $\mathcal{L}$ is the graph Laplacian matrix with $\mathcal{L}_{ij} = -w_{ij} + \delta_{ij}\sum_{k=1}^N w_{ik}$, $I_N$ is the $N\times N$ identity matrix and $\otimes$ is the Kronecker (tensor) product. The eigenvalues of \eqref{eq:Jconstrows} depend on the particular choice of $W$ with constant row sums. In the case of global coupling where $w_{ij}=1/N$ so that $c=1$ we see that $\mathcal{J}$ has non-zero eigenvalues $\kappa +\epsilon (H_5(0) + H_6(0))$ and $\mu_\pm$ each of multiplicity $N-1$ where $\mu_\pm$ are the eigenvalues of
\begin{align}\label{eq:M}
\mathcal{M}(\Psi)&=
\begin{bmatrix}
-\epsilon(H_1'(0) + {\Psi}(H_2'(0) + H_3'(0)) & \epsilon H_2(0) \\
-\epsilon(H_4'(0) + {\Psi}(H_5'(0) + H_6'(0)) & \kappa +\epsilon H_5(0)
\end{bmatrix}.
\end{align}

Therefore synchrony is stable when $\kappa +\epsilon (H_5(0) + H_6(0))<0$, $\mbox{Trace}(\mathcal{M}(\Psi))<0$ and $\det(\mathcal{M}(\Psi))>0$. Reducing back to the phase-only description by taking $H_2, \ldots, H_6 \equiv 0$ we recover the result for phase oscillators that synchrony is stable if $-\epsilon H'_1(0)<0$ \cite{Ashwin1992, Ashwin2016, Coombes2023}.

\paragraph{Interaction functions with $H_1(0)=0$ and $H_4(0)=0$}

If $H_1(0)=0$ and $H_4(0)=0$ then $\Psi= 0$ and therefore the synchronous orbit coincides with the stable uncoupled node orbit and $\Omega = \omega$. The Jacobian is given by
\begin{align}\begin{split}
\mathcal{J} =& -\epsilon\begin{bmatrix}
 H_1'(0) &  H_3(0) \\
 H_4'(0) &  H_6(0)
\end{bmatrix} \otimes \mathcal{L} + \begin{bmatrix}
0 & 0 \\
0 & \kappa
\end{bmatrix} \otimes I_N \\ &+ \begin{bmatrix}
0 & \epsilon (H_2(0)+ H_3(0)) \\
0 & \epsilon  (H_5(0) + H_6(0))
\end{bmatrix} \otimes \mbox{Diag}(c_1, \ldots, c_N).
\end{split}
\end{align} For the case of diffusive coupling which is linear in $\vec{x}_j - \vec{x}_i$, $H_1(0)= H_4(0)$ and additionally, since $J_1= -J_2$, $H_2(0)= - H_3(0)$ and $H_5(0)= - H_6(0)$. This gives
\begin{align}\label{eq:stabsyncpi}
\mathcal{J} =-\epsilon \begin{bmatrix}
 H_1'(0) &  H_3(0) \\
 H_4'(0) &  H_6(0)
\end{bmatrix} \otimes \mathcal{L} + \begin{bmatrix}
0 & 0 \\
0 & \kappa
\end{bmatrix} \otimes I_N,
\end{align} however there are no general results concerning the eigenvalues of $\mathcal{J}$. In the case where coupling is global so that $\mathcal{L}_{ij} = -1/N + \delta_{ij}$ the non-zero eigenvalues of $\mathcal{J}$ are $\kappa$ and $\mu_\pm$ each of multiplicity $N-1$ where $\mu_\pm$ are the eigenvalues of $\mathcal{M}(0)$ in \eqref{eq:M} so that synchrony is stable for global diffusive coupling when
$\kappa + \epsilon H_5(0) -\epsilon H_1'(0)<0$ and $-\epsilon H_1'(0)(\kappa +\epsilon H_5(0)) + \epsilon^2 H_2(0)H_4'(0)>0$.

\subsubsection{The splay state}\label{sec:splay}

Another important example of a $1:1$ phase-locked solution is the splay state, sometimes referred to as a rotating wave or the uniform incoherent state (UIS). Here the phases are uniformly distributed around the unit circle and without loss of generality we take $\theta_{i+1}-\theta_{i}= 2\pi/N$ and $\Phi = (\phi_1, \ldots, \phi_N)$ where $\phi_i=2\pi i/N$. The equations for the collective frequency $\Omega$ and collective orbit $\Psi$, given by \eqref{eqs:phaselocked1} must be satisfied for all values of $i$ for the existence of the splay state. In the case of global coupling, $w_{ij}=1/N$, we find that the splay state exists with 
\begin{align} \label{eq:finiteNsplay}
\Omega = \omega + \frac{\epsilon}{N}\left( \beta_1 + {\Psi}(\beta_2 +\beta_3)\right), \quad
{\Psi} =  -\frac{\epsilon \beta_4}{N\kappa + \epsilon(\beta_5 + \beta_6)},
\end{align} where $\beta_k = \sum_{j=1}^N H_k(\phi_j)$. Considering now the stability of the state we see that the Jacobian takes the form \eqref{eq:pljac} where
\begin{align}
\mathcal{H}_{ij}^{(p)}(\Phi) = \frac{\epsilon}{N} \left( A_{j-i}^{(p)} - \delta_{ij} B^{(p)}\right), \qquad p=1, 2, 3, 4,
\end{align} for
\begin{subequations}
\begin{align}
A_{m}^{(1)} = & H_1'(\phi_m) + {\Psi}\left( H_2'(\phi_m)+ H_3'(\phi_m)\right),\\
B^{(1)} = & \beta'_1 + {\Psi}\left( \beta'_2 +\beta'_3\right), \\
A_{m}^{(2)} = & H_3(\phi_m),\\
B^{(2)} = & -\beta_2, \\
A_{m}^{(3)} = & H_4'(\phi_m) + {\Psi}\left(H_5'(\phi_m)+ H_6'(\phi_m)\right), \\
B^{(3)} =&  \beta'_4 + {\Psi}\left(\beta'_5 + \beta'_6\right),\\
A_{m}^{(4)} = & H_6(\phi_m), \\
B^{(4)} = & -\beta_5 -N\kappa/\epsilon.
\end{align}
\end{subequations}
and $\beta'_k = \sum_{j=1}^N H_k'(\phi_j)$ where $\phi_m = 2\pi m/N$. Therefore the eigenvalues of $\mathcal{H}^{(p)}$ are
\begin{align}
\lambda_q^{(p)} = \frac{\epsilon}{N}\left(\nu_q^{(p)} - B^{(p)}\right), \qquad q=0,1, \ldots, N-1,
\end{align} where $\nu_q^{(p)}$ are the eigenvalues of the circulant matrix $A^{(p)}$ with elements $A_{ij}^{(p)}= A_{j-i}^{(p)}$. Since the $A^{(p)}$ are circulant, the eigenvalues are
\begin{align}
\nu_q^{(p)} = \sum_{m=1}^N A_m^{(p)} \exp(2\pi {\rm i} m q/N),
\end{align} and hence
\begin{subequations}\label{eqs:lambdaq}
\begin{align}\label{eq:lambdaq1}\begin{split}
\lambda_q^{(1)}  = & \frac{\epsilon}{N} \sum_{m=1}^N \left[H_1'(\phi_m) + {\Psi}\left( H_2'(\phi_m)+ H_3'(\phi_m)\right) \right]\\ & \qquad \qquad \times \left( \exp(2\pi {\rm i} m q/N)-1\right),\end{split}\\
\lambda_q^{(2)}  = & \frac{\epsilon}{N} \sum_{m=1}^N \left[ H_3(\phi_m)\exp(2\pi {\rm i} m q/N) + H_2(\phi_m) \right],\\ \label{eq:lambdaq3} \begin{split}
\lambda_q^{(3)}  = & \frac{\epsilon}{N} \sum_{m=1}^N \left[H_4'(\phi_m) + {\Psi}\left( H_5'(\phi_m)+ H_6'(\phi_m)\right) \right]\\ & \qquad\qquad  \times \left( \exp(2\pi {\rm i} m q/N)-1\right),\end{split}\\
\lambda_q^{(4)}  = & \frac{\epsilon}{N} \sum_{m=1}^N \left[ H_6(\phi_m)\exp(2\pi {\rm i} m q/N) + H_5(\phi_m) \right] + \kappa.
\end{align}     
\end{subequations}
Having found the eigenvalues of each of the blocks $\mathcal{H}^{(p)}$ it remains to determine the eigenvalues $\mu$ of the Jacobian \eqref{eq:pljac}. Since each block $\mathcal{H}^{(p)}$ is circulant, they can each be diagonalised as $\mathcal{H}^{(p)} = Q \Lambda^{(p)} Q^{-1}$ where $\Lambda^{(p)} = \mbox{Diag}(\lambda_0^{(p)}, \ldots, \lambda_{N-1}^{(p)})$, and therefore they commute. By \cite[Theorem 1]{Silvester2000}, the eigenvalues $\mu$ then satisfy
\begin{align}
0 & = |\mathcal{J}-\mu I_{2N}| \notag\\
& =  \left| (\mathcal{H}^{(1)} - \mu I_N)(\mathcal{H}^{(4)} - \mu I_N) - \mathcal{H}^{(2)} \mathcal{H}^{(3)}\right|\notag\\
& = |Q|\left| (\Lambda^{(1)} - \mu I_N)(\Lambda^{(4)} - \mu I_N) - \Lambda^{(2)} \Lambda^{(3)}\right||Q^{-1}|\notag\\
&= \left| (\Lambda^{(1)} - \mu I_N)(\Lambda^{(4)} - \mu I_N) - \Lambda^{(2)} \Lambda^{(3)}\right|.
\end{align}    
The characteristic equation for the $2N$ eigenvalues of $\mathcal{J}$ is therefore
\begin{align}
\prod_{q=0}^{N-1}\Bigl[ (\lambda_q^{(1)} - \mu)(\lambda_q^{(4)} - \mu) - \lambda_q^{(2)}\lambda_q^{(3)} \Bigr]= 0.
\end{align} Hence the eigenvalues are $\mu_q^\pm$ for $q=0, \ldots, N-1$ where $\mu_q^\pm$ are the eigenvalues of
\begin{align}\label{eq:splayeigs}\Lambda_q =
\begin{bmatrix} \lambda_q^{(1)} & \lambda_q^{(2)}  \\ \lambda_q^{(3)} & \lambda_q^{(4)}\end{bmatrix}.\end{align} 
The splay state is stable if $\mbox{Re}(\mu_q^\pm)<0$ for all  $q=0, \ldots, N-1$ with the exception of the zero eigenvalue corresponding to perturbations around the limit cycle.

\paragraph{The splay state in the large $N$ limit}
In the thermodynamic limit $N \to \infty$ for global coupling, $w_{ij}=1/N$, the network averages are effectively Riemann sums and so may be replaced by time averages:
\begin{align}
\lim_{N\to \infty} \frac{1}{N} \sum_{j=1}^N H_k(2\pi j/N) = \frac{1}{2\pi} \int_0^{2\pi} H_k(t) \ {\rm d}t := (H_k)_{0}
\end{align} where $H_k$ has a Fourier series representation $H_k(t) = \sum_{n} (H_k)_n \exp({\rm i} n t)$ with $(H_k)_n = \frac{1}{2\pi} \int_0^{2\pi} H_k(t) \exp(-{\rm i}nt){\rm d} t$. Consequently $\Omega$ and ${\Psi}$ from \eqref{eq:finiteNsplay} become
\begin{subequations}
\begin{align} \label{eq:largeNsplay}
\Omega &= \omega + \epsilon\left( (H_1)_0 + {\Psi}((H_2)_0 + (H_3)_0)\right), \\ 
\Psi &=  -\frac{\epsilon (H_4)_0}{\kappa + \epsilon( (H_5)_0 + (H_6)_0)},
\end{align}\end{subequations} and
\begin{subequations}\label{eqs:lambda}
\begin{align}\label{eq:lambda1}
\lambda_q^{(1)}  = & -\epsilon {\rm i} q\left( (H_1)_{-q} + {\Psi}\left( (H_2)_{-q} + (H_3)_{-q} \right) \right),\\
\lambda_q^{(2)}  = & \epsilon\left( (H_3)_{-q} + (H_2)_0 \right),\\
\lambda_q^{(3)}  = & -\epsilon {\rm i} q\left( (H_4)_{-q} + {\Psi}\left( (H_5)_{-q} + (H_6)_{-q} \right) \right),\\
\lambda_q^{(4)}  = &  \kappa + \epsilon\left( (H_6)_{-q} + (H_5)_0 \right)\label{eq:lambda4}.
\end{align}    
\end{subequations} The eigenvalues of the Jacobian $\mathcal{J}$ for the splay state in the large $N$ limit are given by the eigenvalues of \eqref{eq:splayeigs} with $\lambda_q^{(p)}$, $p=1,2,3,4$ now given by \eqref{eqs:lambda}. The splay state is stable if $\mbox{Re}(\mu_q^\pm)<0$ for all non-zero $q$. It is straightforward to show that reducing back to phase-only by letting ${\Psi}=0$ and taking $|\epsilon|\ll |\kappa|$ we recover the result that for globally coupled phase oscillators in the thermodynamic limit, the splay state is stable if $\epsilon q \mbox{Im}((H_1)_{-q})<0$ for all non-zero $q$ \cite{Ashwin1992, Watanabe1997, Ashwin2016, Coombes2023}.

\subsubsection{Balanced cluster states in globally coupled networks}\label{sec:balancedclusters}
Consider now the more general states sharing a common orbit in a globally coupled network, $w_{ij}=1/N$, where the $N=Mm$ nodes form $M$ clusters, each with $m$ nodes. Symmetry imposes that the phase difference between any of the $M$ clusters is an integer multiple of $2\pi/M$ and without loss of generality we order the nodes such that $\Phi = 2\pi(1,\ldots, 1, 2, \ldots, 2, \ldots, M-1, \ldots, M-1, M, \ldots, M)/M$ where each distinct integer value occurs $m$ times. From \eqref{eq:existence} we find that the collective frequency and isostable coordinate of the collective orbit are 
 \begin{align}
        \Omega = \omega + \frac{\epsilon}{M}\left( \sigma_1 + \Psi(\sigma_2 + \sigma_3) \right), \ \ 
        \Psi  = -\frac{\epsilon \sigma_4}{M\kappa + \epsilon(\sigma_5 + \sigma_6)},
    \end{align}
where $\sigma_k = \sum_{j=1}^M H_k(2\pi j/M)$. Considering stability of such a state, the Jacobian is given by \eqref{eq:pljac}, where each block $\mathcal{H}^{(p)}(\Phi)$ has block circulant structure:
\begin{align}\begin{split}
    &\mathcal{H}^{(p)}(\Phi) = \\&\epsilon \begin{bmatrix}
        b_0^{(p)} \mathbbm{1}_{m} - a^{(p)} I_m & b_1^{(p)}\mathbbm{1}_m & \cdots & b_{M-1}^{(p)}\mathbbm{1}_m\\
        b_{M-1}^{(p)}\mathbbm{1}_m & b_0^{(p)} \mathbbm{1}_{m} - a^{(p)} I_m & \cdots & b_{M-2}^{(p)}\mathbbm{1}_m\\
        \vdots & \vdots & \ddots & \vdots\\
        b_1^{(p)}\mathbbm{1}_m & b_2^{(p)}\mathbbm{1}_m & \cdots & b_0^{(p)} \mathbbm{1}_{m} - a^{(p)} I_m
    \end{bmatrix},\end{split}
\end{align}
where $\mathbbm{1}_{m}$ denotes the $m\times m$ matrix of ones and 
\begin{subequations}
\begin{align}
 a^{(1)} &= \frac{1}{M} (\sigma'_1 + \Psi(\sigma'_2 + \sigma'_3))\\
 a^{(2)} &= -\frac{1}{M} \sigma_2 \\
 a^{(3)} &= \frac{1}{M} (\sigma'_4 + \Psi(\sigma'_5 + \sigma'_6))\\
 a^{(4)} &= -\frac{1}{M} \sigma_5 - \kappa/\epsilon  \\
 b_i^{(1)} &= \frac{1}{N} \left(H_1'(2\pi i/M)  + \Psi(H_2'(2\pi i/M)  + H_3'(2\pi i/M))\right)\\
 b_i^{(2)} &= \frac{1}{N} H_3(2\pi i/M) \\
 b_i^{(3)} &= \frac{1}{N} \left(H_4'(2\pi i/M)  + \Psi(H_5'(2\pi i/M)  + H_5'(2\pi i/M))\right)\\
 b_i^{(4)} &= \frac{1}{N} H_6(2\pi i/M)  
\end{align}
\end{subequations}
for $\sigma'_k = \sum_{j=1}^M H_k'(2\pi j/M)$. Due to the block circulant structure, 
\begin{align}
    \det(\mathcal{H}^{(p)}(\Phi)) = \prod_{q=0}^{M-1} \det(\mathcal{J}_q),
\end{align}
where \cite{Tee2007, Kra2012, Coombes2023}
\begin{align}
    \mathcal{J}_q^{(p)} = -\epsilon \left[ a^{(p)}I_m - \sum_{k=0}^{M-1} b_k^{(p)} {\rm e}^{2\pi i k q/M} \mathbbm{1}_m \right].
\end{align} By an application of the matrix determinant lemma \cite{Ding2007} the eigenvalues $\lambda^{(p)}$ and $\lambda^{(p)}_q$ of $\mathcal{H}^{(p)}$ satisfy
\begin{align*}\begin{split}
   0&=  |\lambda^{(p)} I_N - \mathcal{H}^{(p)}|  \\
    & = \prod_{q=0}^{M-1}  \left| (\lambda^{(p)} + \epsilon a^{(p)}) I_m - \epsilon \sum_{k=0}^{M-1} b_k^{(p)} {\rm e}^{2\pi {\rm i} k q/M} \mathbbm{1}_m \right|\\
    &= (\lambda^{(p)} + \epsilon a^{(p)})^{N-M} \\ & \qquad \times \prod_{q=0}^{M-1}  \left( \lambda_q^{(p)} + \epsilon a^{(p)} - \epsilon m \sum_{k=0}^{M-1} b_k^{(p)} {\rm e}^{2\pi {\rm i} k q/M}  \right).
    \end{split}
    \end{align*}
Furthermore, each $\mathcal{H}^{(p)}$ can be diagonalised as $\mathcal{H}^{(p)} = Q \Lambda^{(p)}Q^{-1}$ where 
\begin{align*}
    \Lambda^{(p)} = \mbox{Diag}( \lambda^{(p)}, \ldots, \lambda^{(p)}, \lambda^{(p)}_0, \ldots, \lambda^{(p)}_{M-1} )
\end{align*} with $\lambda^{(p)}$ repeated $N-M$ times. As for the splay state, the $\mathcal{H}^{(p)}$ therefore commute and 
\begin{align*}
0 & = |\mathcal{J}-\mu I_{2N}| \\
&= \left| (\Lambda^{(1)} - \mu I_N)(\Lambda^{(4)} - \mu I_N) - \Lambda^{(2)} \Lambda^{(3)}\right|
\end{align*} so that the characteristic equation of $\mathcal{J}$ is 
\begin{align}
    \begin{split}
        0 &= \left[ (\lambda^{(1)} - \mu)(\lambda^{(4)} - \mu) - \lambda^{(2)}\lambda^{(3)}\right]^{N-M} \\ & \quad \times \prod_{q=0}^{M-1}\left[ (\lambda_q^{(1)} - \mu)(\lambda_q^{(4)} - \mu) - \lambda_q^{(2)}\lambda_q^{(3)}\right]. 
    \end{split}
\end{align} Hence the eigenvalues of $\mathcal{J}$ are the eigenvalues $\mu^\pm$ of 
\begin{align}  \Lambda =
    \begin{bmatrix}
        \lambda^{(1)} & \lambda^{(2)}\\ \lambda^{(3)} & \lambda^{(4)}
    \end{bmatrix} = - \epsilon \begin{bmatrix}
        a^{(1)} & a^{(2)}\\ a^{(3)} & a^{(4)}
    \end{bmatrix}
\end{align} each of multiplicity $N-M$ (corresponding to intracluster perturbations) and the eigenvalues $\mu^\pm_q$ of 
$\Lambda_q$ as in \eqref{eq:splayeigs} for $q=0, \ldots M-1$ where $\lambda_q^{(p)}$ are as in \eqref{eqs:lambdaq} replacing $N$ with $M$ (the intercluster eigenvalues, which are the eigenvalues of the splay state in an $M$ node network). Note that $\Lambda_0$ has a zero eigenvalue (the purely rotational eigenvalue). It is clear that the splay state is a special case where $M=N$, $m=1$. Synchrony is the special case where $M=1$, $m=N$ and we note that the stability conditions for synchrony in globally coupled networks are recovered since here the non-zero eigenvalue of $\Lambda_0$ is $\kappa + \epsilon(H_5(0) + H_6(0))$ and $\Lambda = M(\Psi)$ as in \eqref{eq:M}.

\subsection{Two-cluster states in globally coupled networks}\label{sec:twoclusters}
We now consider network states where the $N$ nodes in a globally coupled network ($w_{ij}=1/N$) lie in two clusters denoted $\mathcal{C}_A$ and $\mathcal{C}_B$ containing $N_A\leq N/2$ and $N_B= N-N_A$ nodes respectively. We label nodes such that $i \in \mathcal{C}_A$ for $i=1, \ldots, N_A$ and $i \in \mathcal{C}_B$ for $i=N_A +1, \ldots, N$. Without loss of generality we may assume that phase of nodes in $\mathcal{C}_A$ is $\theta_A= \Omega t$ and nodes in $\mathcal{C}_B$ have phase $\theta_B = \Omega t + \chi$ where $\Omega$ is the collective frequency of the solution and $\chi$ is the phase difference between the clusters. The clusters will generally have different orbits (coinciding with differing isostables of the uncoupled nodes) such that for $i \in \mathcal{C}_A$,  $\psi_i(t) = \Psi_A$ and for $i \in \mathcal{C}_B$,  $\psi_i(t) = \Psi_B$ where $\Psi_A$ and $\Psi_B$ are constants. 

Using \eqref{eqs:phaselocked1}, equations for the collective frequency, phase difference and cluster orbits can be determined. The four equations defining $\Omega$, $\chi$, ${\Psi}_A$ and ${\Psi}_B$ are \begin{subequations}\label{eqs:clusterexistence}
    \begin{align}
    \begin{split}
\Omega &= \omega + \frac{\epsilon}{N} \Bigl[ N_A\left (H_1(0) + {\Psi}_A(H_2(0) + H_3(0))\right) \\ & \qquad  + N_B\left (H_1(\chi) + {\Psi}_A H_2(\chi) + {\Psi}_B H_3(\chi)\right)\Bigr],\label{eq:thetaA}\end{split}\\ \begin{split}
0 &= \kappa {\Psi}_A + \frac{\epsilon}{N} \Bigl[ N_A\left (H_4(0) + {\Psi}_A(H_5(0) + H_6(0))\right) \\ & \qquad  + N_B\left (H_4(\chi) + {\Psi}_A H_5(\chi) + {\Psi}_B H_6(\chi)\right)\Bigr],\label{eq:psiA}\end{split}\\ \begin{split}
\Omega &= \omega + \frac{\epsilon}{N} \Bigl[ N_A\left (H_1(-\chi) + {\Psi}_B H_2(-\chi) + {\Psi}_A H_3(-\chi) \right) \\ & \qquad + N_B\left (H_1(0) + {\Psi}_B ( H_2(0) + H_3(0)) \right)\Bigr],\label{eq:thetaB} \end{split} \\ \begin{split}
0 &= \kappa {\Psi}_B + \frac{\epsilon}{N} \Bigl[ N_A\left (H_4(-\chi) + {\Psi}_B H_5(-\chi) + {\Psi}_A H_6(-\chi) \right) \\ & \qquad  +  N_B\left (H_4(0) + {\Psi}_B ( H_5(0) + H_6(0)) \right)\Bigr].\label{eq:psiB}
\end{split}\end{align}
\end{subequations}
Equations \eqref{eq:psiA} and \eqref{eq:psiB} can be solved for ${\Psi}_A$ and ${\Psi}_B$ in terms of functions of $\chi$. Substituting the result into \eqref{eq:thetaA} and \eqref{eq:thetaB} we see that $\chi$ is a root of a nonlinear periodic function in terms of functions of $\chi$. The roots can be determined, typically using numerical root finding schemes. For given values of $N_A$, $N_B$ and $\epsilon$ there can be many possible solutions $(\Omega, \chi, {\Psi}_A, {\Psi}_B)$. Synchrony, $(\Omega, 0,{\Psi},{\Psi})$, is a solution for all values of $N_A$, $N_B$ and $\epsilon$ where $\Psi$ and $\Omega$ are given by \eqref{eq:synciso} and \eqref{eq:syncomega} respectively with $c_i =1$ for all $i$. However we are interested in states with $\chi \neq 0$.

Linearising around a a two-cluster state we find that the Jacobian $\mathcal{J}$ is of the form \eqref{eq:pljac} where $\Phi = (0, \ldots, 0, \chi, \ldots \chi)$. Here each of the $N\times N$ block matrices $\mathcal{H}^{(p)}(\Phi)$, $p=1, \ldots, 4$ themselves have a block structure 
\begin{align}\label{eq:2clusterblocks}
\mathcal{H}^{(p)}(\Phi) = 
\left[\begin{array}{@{}c|c@{}}
  M^{(p)}_{AA} & M^{(p)}_{AB} \\
\hline
  M^{(p)}_{BA} & M^{(p)}_{BB}
\end{array}\right],
\end{align} where, letting $\mathbbm{1}_{N\times M}$ denote the $N\times M$ matrix of ones,
\begin{subequations}
\begin{align} 
M^{(p)}_{AA} &= a_p \mathbbm{1}_{N_A \times N_A} + s^A_p I_{N_A}\label{eq:blockAA}\\
M^{(p)}_{AB} &= b_p \mathbbm{1}_{N_A \times N_B}\\
M^{(p)}_{BA} &= c_p \mathbbm{1}_{N_B \times N_A}\\
M^{(p)}_{AA} &= d_p \mathbbm{1}_{N_B \times N_B} + s^B_p I_{N_B}\label{eq:blockBB}.
\end{align}   
\end{subequations} The coefficients $a_p, b_p, c_p, d_p, s^A_p$ and $s^B_p$ are given in Appendix \ref{sec:AppendixB} where they are arranged into $2 \times 2$ matrices $\mathcal{A}$, $\mathcal{B}$, $\mathcal{C}$, $\mathcal{D}$, $\mathcal{S}^A$ and $\mathcal{S}^B$. The block structure of $\mathcal{J}$ facilitates the computation of its eigenvalues. As shown in Appendix \ref{sec:AppendixB} the eigenvalues corresponding to intracluster perturbations are found to be the eigenvalues, $\mu_A^\pm$, of $\mathcal{S}^A$ each with multiplicity $N_A-1$ and the eigenvalues $\mu_B^\pm$ of $\mathcal{S}^B$ each with multiplicity $N_B-1$. For stability with respect to intracluster perturbations it is required that $\mbox{Trace}(\mathcal{S}^A)<0$, $\det(\mathcal{S}^A)>0$, $\mbox{Trace}(\mathcal{S}^B)<0$ and $\det(\mathcal{S}^B)>0$. The remaining $4$ eigenvalues correspond to intercluster perturbations and are the eigenvalues of the $4\times 4$ block matrix 
\begin{align}\label{eq:Jtwoclust}
\mathcal{J}_M =
\begin{bmatrix} \mathcal{S}^A + N_A \mathcal{A} &  N_A\mathcal{B}\\ N_B \mathcal{C} & \mathcal{S}^B + N_B \mathcal{D}\end{bmatrix}.
\end{align} Row and column operations show that one of the eigenvalues of $\mathcal{J}_M$ is the expected zero eigenvalue due to purely rotational eigenmodes. The others, corresponding to changes in the phase difference between clusters, and changes in the two isostable values ${\Psi}_A$, ${\Psi}_B$, are required to have negative real part for stability of the cluster state. They satisfy a cubic characteristic polynomial $\mu^3 + q_1 \mu^2 + q_2 \mu^3 + q_3 =0$. By Routh's stability criterion, all roots of this cubic polynomial lie in the left half-plane when $q_1>0$, $q_2>0$, $q_3>0$ and $q_1 q_2 >q_3$. Reducing to the phase-only case, we recover known results for phase-reduced networks\cite{Kim2000, Brown2003, Coombes2023}.

\section{A comparison of higher order phase and phase-isostable reductions}\label{sec:compare}

In the previous section we have used the phase-isostable network equations \eqref{eqs:piaveraged} to determine with greater accuracy conditions for the existence and stability of phase-locked states. Other authors have used an alternative approach based on isostable coordinates to describe network phenomena which cannot be observed using first-order phase reduction. Rather than work with the phase and isostable equations, it is possible to use isostable coordinates to obtain higher order phase reduction. The process to obtain this reduction, given by Park and Wilson in \cite{Park2021}, involves expanding the isostable coordinate as $\psi(t)= \epsilon p^{(1)}(t) + \epsilon^2 p^{(2)}(t) + \cdots$ where $p^{(k)}(t)$ are $O(1)$. A hierarchy of linear differential equations can be found and solved for the $p^{(k)}(t)$ to any order. Substituting back into the phase equation gives the higher order phase reduced equation including terms up to $O(\epsilon^3)$ and after approximation using first-order averaging as
\begin{widetext}
    \begin{align} \label{eq:avhigherorderphase}
\begin{split}
\FD{\theta_i}{t} = \omega \ +  \ & \epsilon \sum_{j=1}^N w_{ij} \overline{H}_1(\theta_j-\theta_i) +  \epsilon^2 \sum_{j,k=1}^N \left[ w_{ij}w_{ik} \overline{H}_2(\theta_j-\theta_i, \theta_k-\theta_i) + w_{ij}w_{jk} \overline{H}_3( \theta_j-\theta_i, \theta_k-\theta_i) \right] \\  + \  & \epsilon^3 \sum_{j,k, l=1}^N \Bigl[ w_{ij}w_{ik}w_{il} \overline{H}_4(\theta_j-\theta_i, \theta_k-\theta_i, \theta_l-\theta_i)  + w_{ij}w_{ik}w_{kl}  \overline{H}_5( \theta_j-\theta_i, \theta_k-\theta_i, \theta_l-\theta_i) \\ & \qquad \qquad + w_{ij}w_{jk}w_{jl}  \overline{H}_6( \theta_j-\theta_i, \theta_k-\theta_i, \theta_l-\theta_i) + w_{ij}w_{jk}w_{kl}  \overline{H}_7(\theta_j-\theta_i, \theta_k-\theta_i, \theta_l-\theta_i)  \\ & \qquad \qquad \qquad+ w_{ij}w_{ik}w_{jl}  \overline{H}_8( \theta_j-\theta_i, \theta_k-\theta_i, \theta_l-\theta_i)\Bigr], 
\end{split}
\end{align}
\end{widetext} where $\overline{H}_1, \ldots, \overline{H}_8$ are defined in Appendix \ref{sec:AppendixC}. Wilson and Ermentrout \cite{Wilson2019c} previously explicitly derived these equations up to order $\epsilon^2$, while Park and Wilson \cite{Park2021} indicate how the equations may be derived to any order, giving explicit equations to order $\epsilon^3$. However, the form of the higher order phase reduction given in \cite{Park2021} differs from \eqref{eq:avhigherorderphase}. In line with other work on higher order phase equations \cite{Leon2019, Battison2020} we observe in \eqref{eq:avhigherorderphase} that for networks of more than two nodes the higher order terms have non-pairwise phase interactions (i.e., terms involving the phases of three or more oscillators), despite the interactions between the nonlinear oscillators in \eqref{eq:fullnetwork} being pairwise. Park and Wilson \cite{Park2021} however arrive at a higher order phase reduced equation with pairwise interaction terms which we find to be erroneous. We give details of our calculations leading to \eqref{eq:avhigherorderphase} in Appendix \ref{sec:AppendixC}. More recently Park and Wilson \cite{Park2023} have updated their calculations allowing for non-identical oscillators. 

We therefore see that there are two related but different frameworks, both using isostable coordinates, which extend standard first-order phase reduction; the phase-isostable network equations \eqref{eqs:piaveraged} and the higher order phase reduction \eqref{eq:avhigherorderphase}. An obvious question is which approach can most accurately capture the dynamics of the full system \eqref{eq:fullnetwork}? In this section we answer this question for the archetypal example of the mean-field complex Ginzburg-Landau equation where the linear stability boundaries for synchrony and the splay state are known and may be compared with approximations of the boundaries found using the phase-isostable network equations \eqref{eqs:piaveraged} and higher order phase equations \eqref{eq:avhigherorderphase}.

\subsection{The mean-field complex Ginzburg-Landau equation}\label{sec:MFCGL}
The normal form of the Hopf bifurcation (or Stuart-Landau oscillator) is a ubiquitous example of an oscillatory system. When globally diffusively coupled the result is the mean-field complex Ginzburg-Landau equation (MF-CGLE). In real coordinates the governing equations are of the form \eqref{eq:fullnetwork} where $\vec{x}_i = (x_i, y_i)^T$, $w_{ij}=1/N$
\begin{align}\begin{split}\label{eq:MFCGL}
\vec{F}(\vec{x}_i) &= \begin{bmatrix} x_i - (x_i - c_2 y_i)(x_i^2 + y_i^2)\\ y_i - (y_i + c_2 x_i)(x_i^2 + y_i^2)
\end{bmatrix} \\  \vec{G}(\vec{x}_i, \vec{x}_j) &=  \begin{bmatrix} x_j - x_i - c_1(y_j-y_i)\\ y_j - y_i + c_1(x_j-x_i)
\end{bmatrix}\end{split}
\end{align} and $c_1$ and $c_2$ are two real-valued parameters. When $\epsilon=0$ each uncoupled node has a stable limit cycle $\vec{x}^\gamma(t) = \left( \cos(c_2 t), -\sin(c_2 t)\right)^T$ with period $T=2\pi/c_2$ and Floquet exponent $\kappa = -2$. We define the phase $\theta_j$ on the limit cycle as $\theta_j = c_2 t + \phi_j$ so that $\theta_j \in [0, 2\pi)$ and $\omega = c_2$.

In the full system \eqref{eq:MFCGL}, the synchronous state, where $\vec{x}_j = \vec{x}^\gamma(t)$ for $j = 1, \ldots, N$, exists for all parameter values. The splay state is given by $\vec{x}_j= (r \cos(\varphi_j), r \sin(\varphi_j))^T$ where $r= \sqrt{1-\epsilon}$ and $\varphi_j = (-c_2 + \epsilon(c_2-c_1))t + 2\pi j/N$ and exists when $\epsilon<1$. The linear stability analysis of these solutions gives closed formulae for the marginal stability of synchrony (denoted $\epsilon_s$) and the splay state (denoted $\epsilon_0$ when $N\geq 3$ and $\epsilon_a$ when $N=2$) for fixed values of $c_1$ and $c_2$ as \cite{Hakim1992, Nakagawa1993}
\begin{subequations}
    \begin{align}\label{eq:syncstab}
\epsilon_s = \frac{-2(1+c_1 c_2)}{1+c_1^2},
\end{align}
\begin{align}\label{eq:antisyncstab}
&\epsilon_a^2 (c_1^2 + 2c_1c_2 +3) - 2 \epsilon_a (1+c_1c_2)=0, \\ \begin{split}  \label{eq:splaystab}
&\epsilon_0(2 \epsilon_0-1)c_1^2 + 4(\epsilon_0-1)(2\epsilon_0-1)c_1c_2 \\ & \qquad - \epsilon_0(\epsilon_0-1)c_2^2  + (3\epsilon_0-2)^2 =0.\end{split} 
\end{align}
\end{subequations} Here we compute and compare the approximations to the marginal stability curves for these states using both the phase-isostable network equations \eqref{eqs:piaveraged} and the higher order phase reduction \eqref{eq:avhigherorderphase} in order to compare their accuracy to the true curves. For completeness the results are also compared with those obtained by Le\'{o}n and Paz\'{o} \cite{Leon2019} who obtained an alternative higher order phase reduction based on the isochrons which are known in closed form for this system \cite{Castejon2013, Monga2019}. We note that both methods based on isostable reduction will generalise to any system, however the approach of \cite{Leon2019} relies on having a closed formula for the isochrons. We choose to compare the phase-isostable network equations to first-order in $\psi$ with the second and third order phase reductions since these reductions require comparable computational effort in terms of calculating the required terms in the expansions \eqref{eqs:expansions}.

 Using \eqref{eqs:adjoints}, normalising as appropriate and denoting $\hat{r}(t) = (\cos(c_2t), -\sin(c_2t))^T$ and $\hat{\varphi}(t) = (\sin(c_2t), \cos(c_2t))^T$ we obtain
\begin{align*}
 g^{(1)}(t) &= A(\hat{r}(t) + c_2 \hat{\varphi}(t)), \\
g^{(2)}(t) &=\frac{A^2}{2} \left( (3-c_2^2) \hat{r}(t) + 4c_2 \hat{\varphi}(t)\right), \end{align*} \begin{align*}
Z^{(0)}(t) &= c_2 \hat{r}(t) - \hat{\varphi}(t),\\ 
Z^{(1)}(t) & = \frac{1}{A}\hat{\varphi}(t),\\
Z^{(2)}(t) & = -\frac{c_2}{2}\hat{r}(t) + \frac{1}{2}\hat{\varphi}(t),\end{align*} \begin{align*}
I^{(0)}(t) &= \frac{1}{A}\hat{r}(t),\\ 
I^{(1)}(t) & = - 3\hat{r}(t) + c_2 \hat{\varphi}(t), \\
I^{(2)}(t) & = \frac{A}{2}\left((3-c_2^2)\hat{r}(t) -4c_2\hat{\varphi}(t)\right),
\end{align*} where $A = \|g^{(1)}(0)\|= (1+c_2^2)^{-1/2}$. Calculating the expansions for the coupling function \eqref{eq:couplingexp}, we find that $J_2=-J_1$ and since the coupling is linear, all Hessian terms $\mathbf{H}_q$ vanish.

\subsubsection{The phase-isostable network equation transformation}
Considering first the transformation to phase-isostable network equations \eqref{eqs:piaveraged}, we calculate that the averaged coupling functions as in \eqref{eq:Hk} are
\begin{subequations}
\begin{align}
H_1(\chi)&= (c_2-c_1)\left(\cos(\chi)-1\right) + (1+c_1c_2)\sin(\chi),\\
H_2(\chi)&= A(1+c_2^2)(c_1\cos(\chi)-\sin(\chi)) = - H_3(\chi),\\
H_4(\chi)&= \frac{1}{A} \left(c_1 \sin(\chi) + \cos(\chi) -1\right),\\
H_5(\chi)&= 2 + (c_1c_2 -3) \cos(\chi) - (3c_1 +c_2) \sin(\chi),\\
H_6(\chi) &= (c_1+c_2) \sin(\chi) + (1-c_1c_2)\cos(\chi).
\end{align}
\end{subequations}
From \eqref{eq:M} with $\Psi=0$, we find that synchrony is stable in the globally coupled network when $\epsilon>-1$ and
\begin{align*}
2\epsilon (1+c_1c_2) + \epsilon^2(1+c_1^2)>0.
\end{align*} The linear stability boundary for synchrony using the phase-isostable network equations, denoted $\epsilon_{s, PI}$ is therefore identical to the closed formula for the full dynamics \eqref{eq:syncstab} $\epsilon_{s, PI} =\epsilon_s$.

For the splay state (that exists only when $\epsilon<1$), we find that the collective orbit satisfies ${\Psi}= \epsilon/(2A(\epsilon-1))$ and $\Omega = c_2 - \epsilon(c_2-c_1)$ from \eqref{eq:finiteNsplay}. It can be seen that this agrees precisely with the orbit in the full system \eqref{eq:MFCGL} as it can be computed that for the MF-CGLE, isostables are circles of radius $r(\psi) = \sqrt{k/(\psi +k)}$ where $I^{(0)} = -2k\hat{r}(t)$ \cite{Castejon2013, Monga2019}. Since we have chosen to normalise such that $k=-(2A)^{-1}$ we find that the orbit for the splay state in the phase-isostable transformed system is  $r({{\Psi}}) = \sqrt{1-\epsilon}$.

Using the results of section \ref{sec:splay} we can compute the conditions for linear stability of the splay state. When $N=2$ (so the splay state is the antisynchronous solution) the eigenvalues of the Jacobian are $0$, $-2(1-\epsilon)$ and the
eigenvalues of $\Lambda_1$ given by \eqref{eq:splayeigs}. Since the matrix $\Lambda_1$ has real entries in this case, the eigenvalues have negative real part when $\mbox{Tr}(\Lambda_1)<0$ and $\mbox{Det}(\Lambda_1)>0$. This gives stability boundaries $\epsilon=1/2$, $\epsilon=0$ and $\epsilon= \epsilon_{a,PI}$ where $\epsilon_{a,PI}$ satisfies
\begin{align}\label{eq:antiPI}\begin{split} 
 (c_1^2 c_2^2 - 2c_1c_2 -3)\epsilon_{a,PI}^2 & + (4 c_1c_2 + c_1^2 +5)\epsilon_{a,PI} \\ & - 2(c_1c_2 +1) =0.\end{split}
\end{align}

When $N\geq 3$ and finite and also in the large $N$ limit the non trivial marginal linear stability boundary for the splay state, denoted $\epsilon= \epsilon_{0, PI}$, can be shown to be given by (see Appendix \ref{sec:AppendixD}), \begin{align}\label{eq:splaystabboundary}
\begin{split}
&\epsilon^5\left( (c_2^2 +9)(1+c_1c_2)(c_1c_2-5)\right) \\ & + \epsilon^4\left(8 c_2^3 c_1 + (5-19c_1^2) c_2^2 + 152 c_1 c_2 +9 c_1^2 +177 \right) \\ & + \epsilon^3 \left(-4c_2^3 c_1  + 8(2c_1^2 +1)c_2^2-260 c_1 c_2 - 20 c_1^2 -284\right) \\ & + \epsilon^2 \left( -4(3+c_1^2)c_2^2 + 224 c_1 c_2 + 16 c_1^2 +232\right)\\ & + \epsilon\left(4(c_2^2-c_1^2) - 96(1+c_1c_2) \right)  + 16(1+c_1 c_2)=0.\end{split}
\end{align}

\subsubsection{The higher-order phase reduction}

The synchronous solution of \eqref{eq:avhigherorderphase} with $\theta_i= \theta$ for $i=1, \ldots, N$ in the case of global coupling $w_{ij}=1/N$ is guaranteed to exist. The linearisation about the synchronous state has the Jacobian
\begin{align}
\begin{split}
-\mathcal{L}\xi  = &- \mathcal{L} \left(\epsilon \overline{H}_1'(0) + \epsilon^2 \sum_{m=2}^3 \left[ \left. \PD{\overline{H}_m}{\chi}\right|_{(0,0)} + \left. \PD{\overline{H}_m}{\eta}\right|_{(0,0)}\right]\right. \\ &\left.  +  \epsilon^3 \sum_{m=4}^8 \left[ \left. \PD{\overline{H}_m}{\chi}\right|_{(0,0,0)} + \left. \PD{\overline{H}_m}{\eta}\right|_{(0,0,0)} +  \left. \PD{\overline{H}_m}{\xi}\right|_{(0,0,0)}\right]\right),\end{split}
\end{align} where $\mathcal{L}$ is the graph Laplacian for global coupling given by $\mathcal{L}_{ij}=-1/N + \delta_{ij}$. This has an ($N-1$ degenerate) eigenvalue of $+1$ and therefore synchrony is stable in the higher-order phase reduction if $\xi>0$. Calculating each of the functions $\overline{H}_1, \ldots, \overline{H}_8$ for the MF-CGLE as described in Appendix \ref{sec:AppendixC} we observe that
\begin{align}\label{eq:hostabsync}
\xi = \epsilon (1+c_1c_2) + \epsilon^2\frac{(1+c_2^2)c_1^2}{2} + \epsilon^3 \frac{(c_1c_2-1)(1+c_2^2)c_1^2}{4}.
\end{align} Therefore, keeping only the terms up to $\epsilon^2$ we obtain the linear stability boundary for the second-order phase reduction as
\begin{align}
\epsilon_{s,2} = \frac{-2(1+c_1c_2)}{c_1^2(1+c_2^2)},
\end{align} which matches identically with the result from the second order phase reduction of Le\'{o}n and Paz\'{o} \cite{Leon2019} whose reduction is based on the isochrons for the MF-CGLE rather than the isostable coordinates used here. The third-order phase reduction gives the values for marginal linear stability for fixed values of $c_1$ and $c_2$ as $\epsilon_{s,3}$ satisfying
\begin{align}
\epsilon_{s,3}^2 (c_1c_2-1)(1+c_2^2)c_1^2 + 2\epsilon_{s,3}(1+c_2^2)c_1^2 + 4(1+c_1c_2)=0,
\end{align} whereas the third-order reduction of Le\'{o}n and Paz\'{o} \cite{Leon2019} gives the value as $\epsilon_{s,3}^*$ satisfying
\begin{align}
(\epsilon_{s,3}^*)^2 (1+c_2^2)c_1^3c_2 + \epsilon_{s,3}^*(1+c_2^2)c_1^2 + 2(1+c_1c_2)=0.
\end{align}

Considering next the stability of the splay state for $N\geq 3$, we find that the Jacobian has the form
\begin{align}
\mathcal{J}_{ij} = A_{j-i} - \delta_{ij} \Upsilon,
\end{align}
where \begin{widetext}\begin{align}\begin{split}
A_j =&  \frac{\epsilon}{N} \overline{H}_1'(2\pi j/N)  + \frac{\epsilon^2}{N^2} \sum_{m=2}^3 \sum_{k=1}^N \left(\PD{\overline{H}_m}{\chi}(2\pi j/N, 2\pi k/N) + \PD{\overline{H}_m}{\eta}(2\pi k/N, 2\pi j/N)\right)\\  & + \frac{\epsilon^3}{N^3} \sum_{m=4}^8 \sum_{k,l=1}^N \left(\PD{\overline{H}_m}{\chi}(2\pi j/N, 2\pi k/N, 2\pi l/N) + \PD{\overline{H}_m}{\eta}(2\pi k/N, 2\pi j/N, 2\pi l/N) + \PD{\overline{H}_m}{\xi}(2\pi k/N, 2\pi l/N, 2\pi j/N)\right), \end{split}
\end{align} \end{widetext} and $\Upsilon = \sum_{k=1}^N A_k$. The circulant structure allows us to determine that the eigenvalues of $\mathcal{J}$ are $\lambda_p= \nu_p-\Upsilon$ where $\nu_p = \sum_{j=1}^N  A_j \exp(2\pi i j p/N)$ for $p=0, \ldots, N-1$. Then
\begin{align}
\lambda_p = \sum_{j=1}^N  A_j \left(\exp(2\pi i j p/N)-1\right), \qquad p=0, \ldots, N-1.
\end{align} Since $\lambda_0=0$ the splay state is stable if $\mbox{Re}(\lambda_p)<0$ for all $p\neq 0$.

When $N\geq 3$, the second-order phase reduction gives the boundary for linear stability of the splay state as
\begin{align}\label{eq:e02}
\epsilon_{0,2} = \frac{4(1+c_1c_2)}{(c_1^2-1)(1+c_2^2)},
\end{align} which again agrees with the second-order phase result of \cite{Leon2019}. To third-order, we find the approximation to the boundary is given by
\begin{align}\label{eq:e03}\begin{split}
\epsilon_{0,3}^2 (c_2^2 +1)(c_2c_1^3 - 3c_1c_2 - 7c_1^2 +5)& \\ + 4 \epsilon_{0,3} (1-c_1^2)(1+c_2^2) + 16(1+c_1 c_2)&=0
\end{split}\end{align} which differs from the result of \cite{Leon2019}, who find the approximation
\begin{align}\begin{split}
(\epsilon_{0,3}^*)^2 (1+c_2^2) (2-2c_1^2 - 3c_1c_2 + c_1^3 c_2)& \\  +  2\epsilon_{0,3}^* (1+c_2^2)(1-c_1^2) + 8(1+c_1c_2)& =0.
\end{split}\end{align}

When $N=2$ we find that for the antisynchronous solution the stability boundaries for the second and third order phase reductions are respectively
\begin{align}
\epsilon_{a,2} = \frac{2(1+c_1c_2)}{c_1^2(1+c_2^2)},\end{align} \begin{align}
\epsilon_{a,3}^2 (1+c_2^2) c_1^2(c_1c_2-3) -  2\epsilon_{0,3} (1+c_2^2)c_1^2+ 4(1+c_1c_2)=0.
\end{align}

\subsubsection{Comparison of approaches}

\begin{figure*}
\begin{center}
\includegraphics[width = \textwidth]{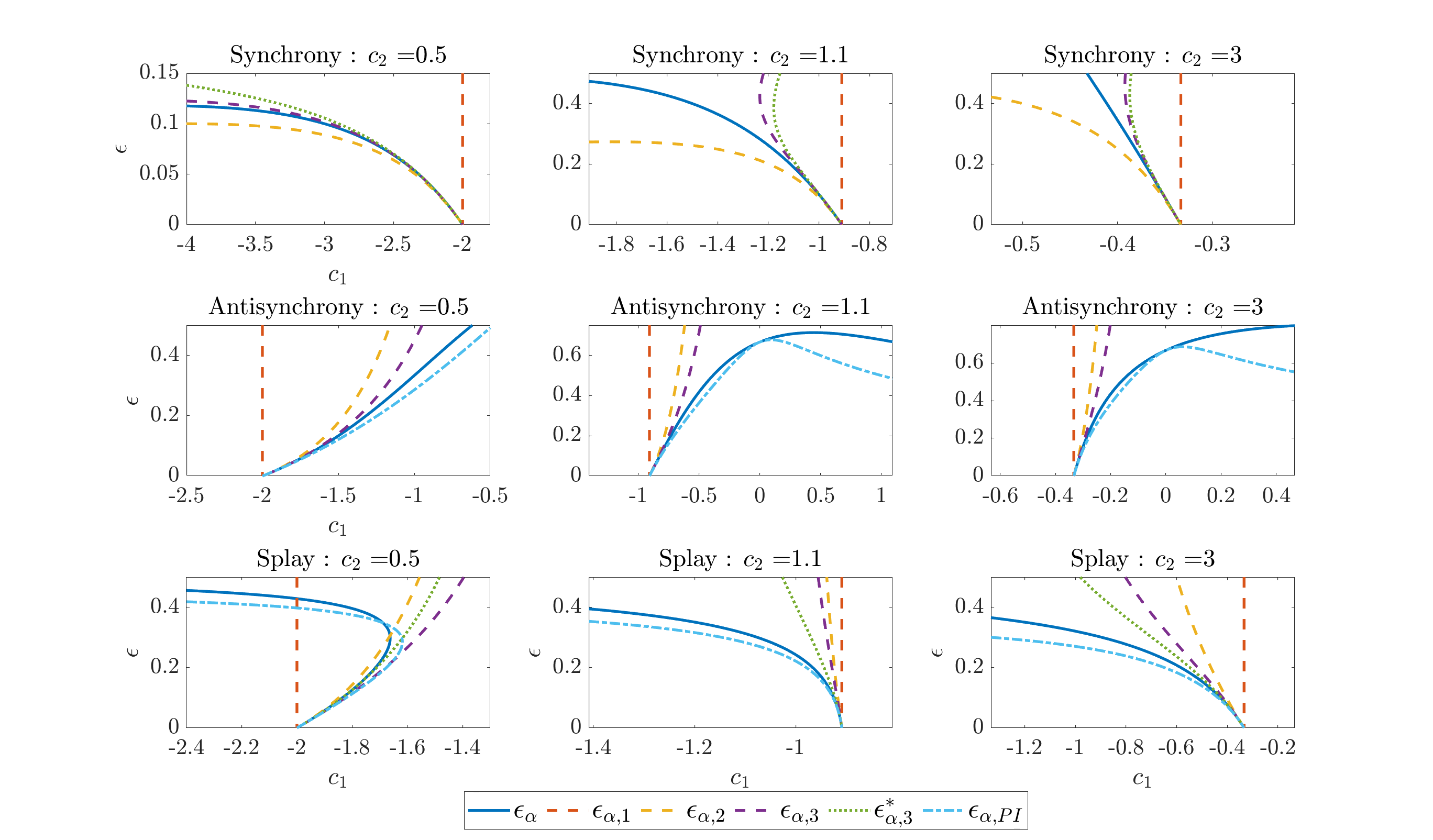}
\caption{Comparison of the marginal stability curves for synchrony, ($\alpha=s$, top), the antisynchronous solution when $N=2$, ($\alpha = a$, middle), and splay state when $N\geq 3$ ($\alpha =0$, bottom) in the MF-CGLE ($\epsilon_{\alpha}$, solid blue); first ($\epsilon_{\alpha, 1}$, dashed red), second ($\epsilon_{\alpha, 2}$, dashed yellow) and third ($\epsilon_{\alpha, 3}$, dashed purple) order phase reductions; third order reduction of \cite{Leon2019} ($\epsilon_{\alpha, 3}^*$, dotted green) and the phase-isostable network approximation ($\epsilon_{\alpha, PI}$, dash-dotted light blue).  For synchrony, the phase-isostable network approximation $\epsilon_{s, PI}$ agrees precisely with the exact MF-CGLE boundary. The second order reduction of \cite{Leon2019} agrees with the second order phase reduction for synchrony and the splay state and $\epsilon_{\alpha, 2}^*$ and $\epsilon_{\alpha, 3}^*$ are not computed for the antisynchronous solution.
Left: $c_2=0.5$, centre: $c_2=1.1$, right: $c_2=3$.
 }
\label{fig:MFCGL}
\end{center}
\end{figure*}

In Figure~\ref{fig:MFCGL} we illustrate for comparison the stability boundary for each of the approaches (exact, phase-isostable network equations, phase interactions of first, second and third order and the approximations of \cite{Leon2019}) of each phase-locked state (synchrony, antisynchrony when $N=2$, splay state for $N\geq 3$) when $\epsilon$ is positive. We choose to fix values of $c_2$ (the intrinsic parameter for the node dynamics determining the angular velocity) and plot the marginal stability curves in the $(c_1, \epsilon)$ plane, following \cite{Nakagawa1993, Leon2019}. Figure \ref{fig:MFCGL} shows the boundaries for $c_2=0.5$ (left), $c_2=1.1$ (centre) and $c_2=3$ (right). Note that we avoid making the choice $c_2=1$ since in this case \eqref{eq:e02} and \eqref{eq:e03} both have a factor of $1+c_1$ giving $c_1=-1$ as a stability boundary for both higher-order phase approximations.

For synchrony we have already observed exact agreement of the stability boundary for the MF-CGLE and the phase-isostable network approximation. The top row of Figure \ref{fig:MFCGL} shows that the third order phase approximation ($\epsilon_{s, 3}$ in dashed purple), provides a slightly better approximation than the third order result of Le\'{o}n and Paz\'{o} ($\epsilon_{s, 3}^*$, in dotted green) \cite{Leon2019}, in a neighbourhood of $(c_1, \epsilon) = (-1/c_2, 0)$.

For the antisynchronous solution when $N=2$ (middle row of Figure \ref{fig:MFCGL}) we see that phase-isostable approximation ($\epsilon_{a,PI}$ of \eqref{eq:antiPI}, shown in dash-dotted light blue) most closely matches the curve for the MF-CGLE over a range of values of $c_1>-1/c_2$ and $\epsilon>0$.

For the splay state where $N\geq 3$ (bottom row of Figure \ref{fig:MFCGL}), while $\epsilon_{0, 3}^*$ (dotted green) may provide the closest approximation to the curve for the MF-CGLE in a neighbourhood of $(c_1, \epsilon) = (-1/c_2, 0)$ we observe that this curve (and those for higher-order phase shown by dotted curves) blow-up for values of $c_1$ away from $-1/c_2$, while the curve for the phase-isostable approximation $\epsilon_{0, PI}$, shown in dash-dotted light blue, is the only one which has a similar shape to the curve for the MF-CGLE, providing a bifurcation locus for all values of $c_1<-1/c_2$. Therefore the phase-isostable approximation provides the most accurate description of the qualitative behaviour (i.e., the existence of a bifurcation) over a range of values of $c_1$ and moderate values of $\epsilon$ although it does not accurately predict the value of $\epsilon$ at which the bifurcation will occur, underestimating in all cases shown in Figure~\ref{fig:MFCGL}.

\begin{figure}
\begin{center}
\includegraphics[width = 0.45\textwidth]{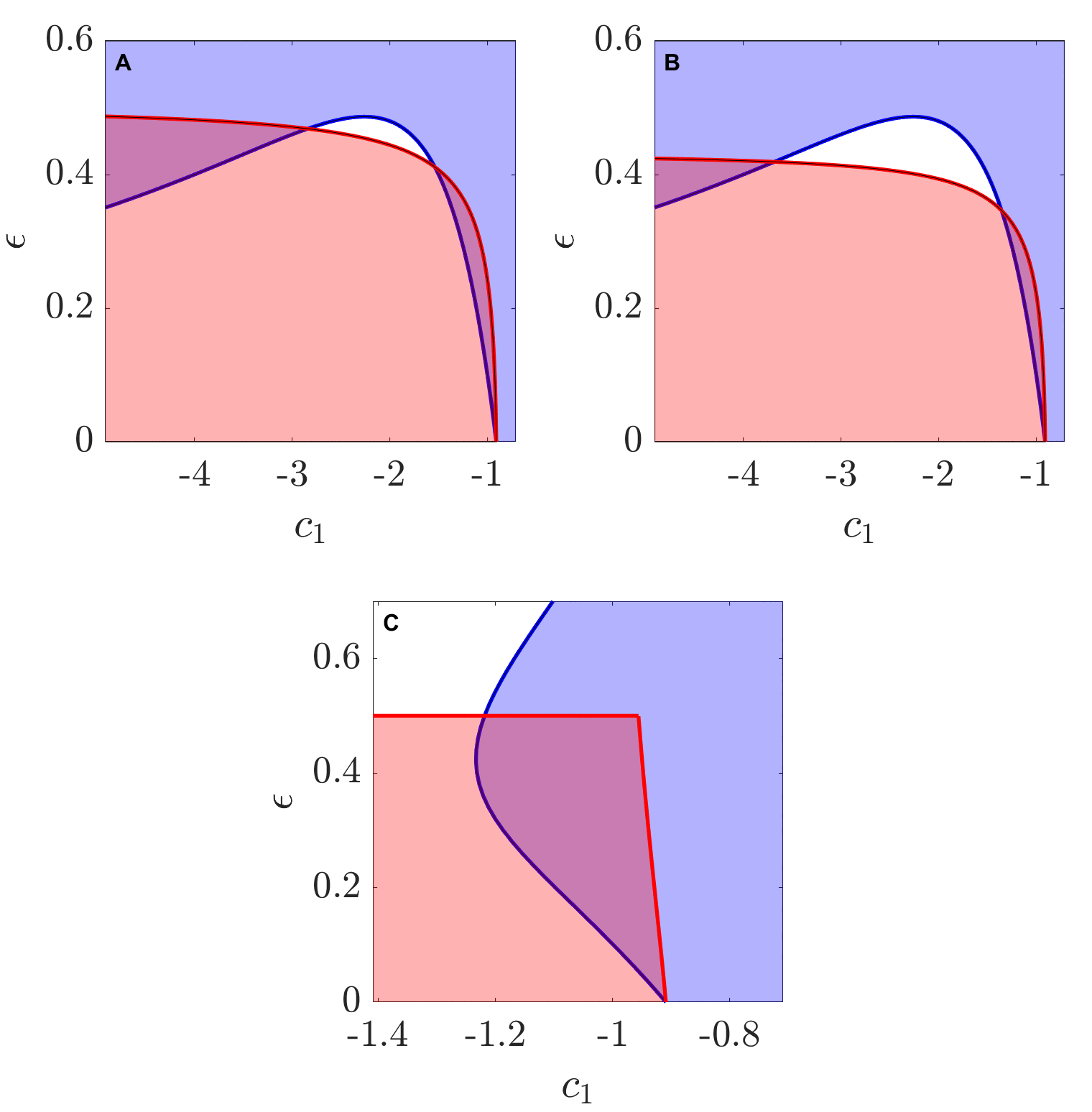}
\caption{Stability regions for synchrony (blue) and splay state (red) for (A) the MF-CGLE, (B) the phase-isostable approximation, (C) the third-order phase approximation when $c_2=1.1$. The phase-isostable network equations (B) can qualitatively reproduce the regions of bistability observed for the full system (A), but higher-order phase approximations cannot (C).
 }
\label{fig:bistability}
\end{center}
\end{figure}

We further note that in the MF-CGLE, below the critical value of $c_2 = \sqrt{3}$ at which the boundaries $\epsilon_s$ and $\epsilon_0$ become tangent at $\epsilon=0$, there are regions of bistability between synchrony and the splay state \cite{Leon2019}, as shown in Figure~\ref{fig:bistability}A for $c_2=1.1$. The phase-isostable approximation can reproduce the bistability regions qualitatively, (see Figure~\ref{fig:bistability}B) while the phase-based approximations and those of \cite{Leon2019} cannot. See Figure~\ref{fig:bistability}C for the case of the third-order phase approximation when $c_2=1.1$ which has very different regions of bistability.

\begin{figure*}
\begin{center}
\includegraphics[width = \textwidth]{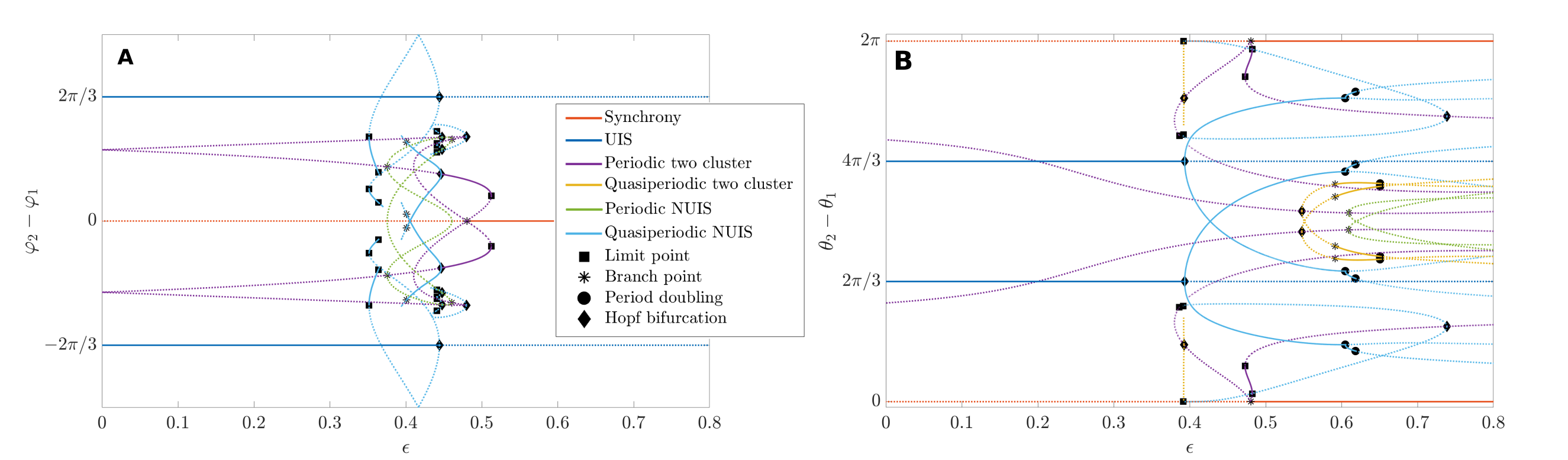}
\caption{Bifurcation diagrams for the MF-CGLE when $N=3$ with interaction strength $\epsilon$ as the bifurcation parameter. (A) shows the result for the full model equations \eqref{eq:MFCGL} with the angular difference between two nodes on the vertical axis. (B) shows the approximation using the phase-isostable network equation framework \eqref{eqs:piaveraged} with the phase difference $\theta_2- \theta_1$ on the vertical axis. 
Solid (dotted) lines indicate stable (unstable) solutions. Shown are two of the three symmetric branches of each type. We do not include the branches with $\theta_2=\theta_1$ since these would obscure the synchronous branch. Parameter values are $c_1=-2$ and $c_2=1.1$. The bifurcation diagrams are produced using the implementation of AUTO in XPPAUT \cite{XPPAUT} and bifurcation locations in (B) are in agreement with stability computations is section \ref{sec:phaselocked}.}
\label{fig:MFCGLbif}
\end{center}
\end{figure*}

Both the full model and the phase-isostable approximation have a parameter region where neither synchrony nor the splay state are stable when $c_2=1.1$. For $c_1=-2$ we compare the stable behaviour in this region for the full model with that predicted by the phase-isostable approximation for a network of $N=3$ nodes. The bifurcation diagrams are as in Figure~\ref{fig:MFCGLbif}. We observe that the phase-isostable approximation captures the bifurcations of synchrony and the splay state (UIS) at moderate values of $\epsilon$ which cannot be seen with a first-order phase reduction. Furthermore, the phase-isostable approximation correctly shows the loss of stability of synchrony to a periodic two-cluster state as $\epsilon$ decreases through $0.48$. It also correctly predicts that as $\epsilon$ increases the splay state losses stability at a Hopf bifurcation to a quasiperiodic nonuniform incoherent state (NUIS) where the three nodes all have different trajectories. Again we observe that the phase-isostable approximation captures the qualitative dynamics, despite discrepancies in the precise values of $\epsilon$ for which some of the bifurcations occur and the precise branching structure. For instance in the full model the Hopf bifurcation of the splay state occurs at $\epsilon=0.393$ which in the full model this occurs at $\epsilon=0.444$. The phase-isostable approximation shows the stable quasiperiodic NUIS bifurcating supercritically directly from the splay state, while in reality the Hopf bifurcation is subcritical and the stable quasiperiodic NUIS bifurcates from a periodic two-cluster state. The phase-isostable approximation also predicts increasingly complex branching of unstable solutions for increasing $\epsilon>0.6$  is not found in the full model. However, for values of $\epsilon$ this large the phase-isostable approximation is not expected to be valid. 

Park and Wilson \cite{Park2023} show analysis of a similar three node network of diffusively coupled complex Ginzburg-Landau models using the higher-order phase reduction. They also find the loss of stability of the splay state at a Hopf bifurcation to a quasiperiodic NUIS with the third-order phase-reduction though they do this through numerical simulation and do not provide any bifurcation diagrams. Our analysis highlights that the interaction functions can be calculated analytically and our computations of stability conditions for phase-locked states in section \ref{sec:phaselocked} allows for explicit location of bifurcations. 

We conclude that the phase-isostable network transformation better captures the qualitative bifurcation structure of the MF-CGLE when compared with second and third order phase reductions. Since it is the normal form for globally diffusively coupled Hopf normal form oscillators, it is natural to suppose that this framework is the most appropriate to use to study linearly coupled networks of oscillatory nodes near Hopf bifurcation. 

It is noted by Park and Wilson \cite{Park2021} that the higher-order phase approximation to the analytical stability boundaries for synchrony and the antisynchronous solution for a network of $N=2$ nodes does improve with increasing order. At 10th order good quantitative and qualitative agreement of the stability boundaries is achieved for moderate values of $\epsilon$, however this does come at huge computational expense in terms of determining all required terms in the PRC and IRC expansions. Using the phase-isostable network equation approximation to first order may be an acceptable compromise, to be able to indicate the qualitative network behaviour using a simpler lower order approximation.

\section{Networks of Morris-Lecar neurons}\label{sec:ML}

\begin{figure}
\begin{center}
\includegraphics[width = 0.49\textwidth]{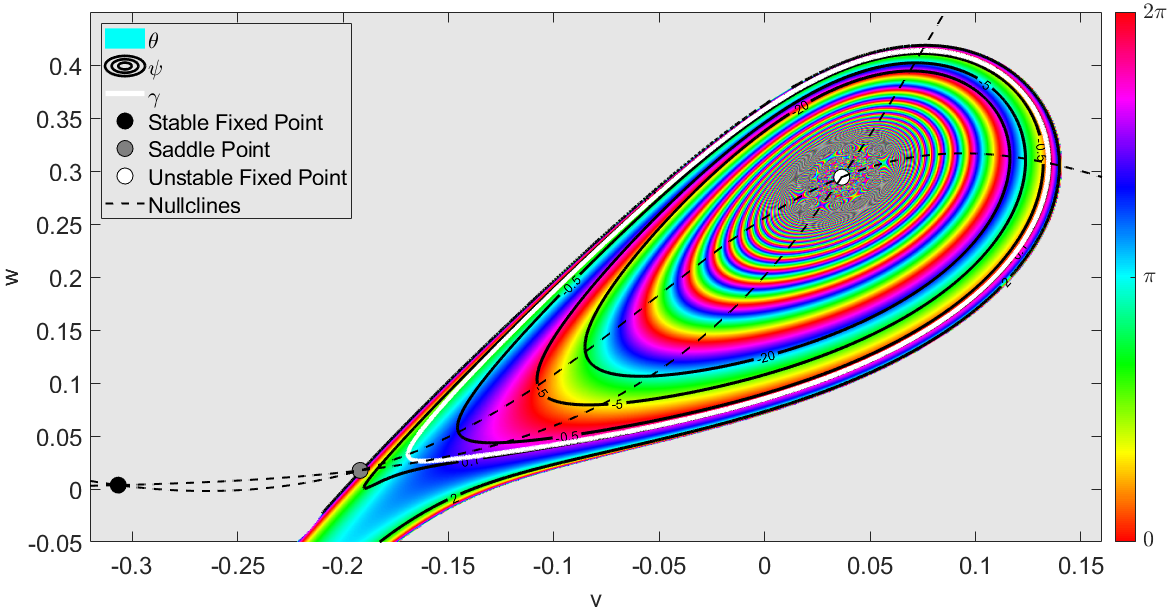}
\caption{Phase portrait for the Morris-Lecar model \eqref{eqs:ML} indicating the stable limit cycle (in white). Parameter values are as given in Appendix \ref{sec:AppendixE}. These parameters place the dynamics near the homoclinic bifurcation so that solutions on the limit cycle spend a significant portion of the period near the saddle point. Also shown are the boundary of the basin of attraction of the limit cycle, a selection of isostables and colours indicate isochrons connecting points with the same asymptotic phase. These are computed numerically using the constructive definition \eqref{eq:isostablesWilson} and the method of Fourier averages \cite{Mauroy2012} respectively. }
\label{fig:MLphaseplane}
\end{center}
\end{figure}

The comparison of higher-order phase reductions and the phase-isostable network equations for the MF-CGLE in section \ref{sec:compare} has shown that the phase-isostable network equation approach is able to predict the qualitative network behaviour for moderate values of interaction strength $\epsilon$, which the higher-order phase reductions of similar computation expense cannot. We now use the phase-isostable framework to reveal network behaviours in small and large networks of planar Morris-Lecar neurons \cite{Morris1981} which cannot be described using traditional first-order phase reduction. The Morris-Lecar model is planar model of neuronal excitability in which oscillations can occur through Hopf, SNIC or homoclinic bifurcations. The equations describing the node dynamics and our choices of parameter values are given in Appendix \ref{sec:AppendixE}. Our parameter choices place the node dynamics near the homoclinic bifurcation. The stable periodic orbit has $\kappa =-0.4094$ and is indicated in Figure \ref{fig:MLphaseplane} along with isochrons and selected isostables normalised so that negative values of the isostable coordinate correspond to points inside the limit cycle. 

\begin{figure}
\begin{center}
\includegraphics[width = 0.35\textwidth]{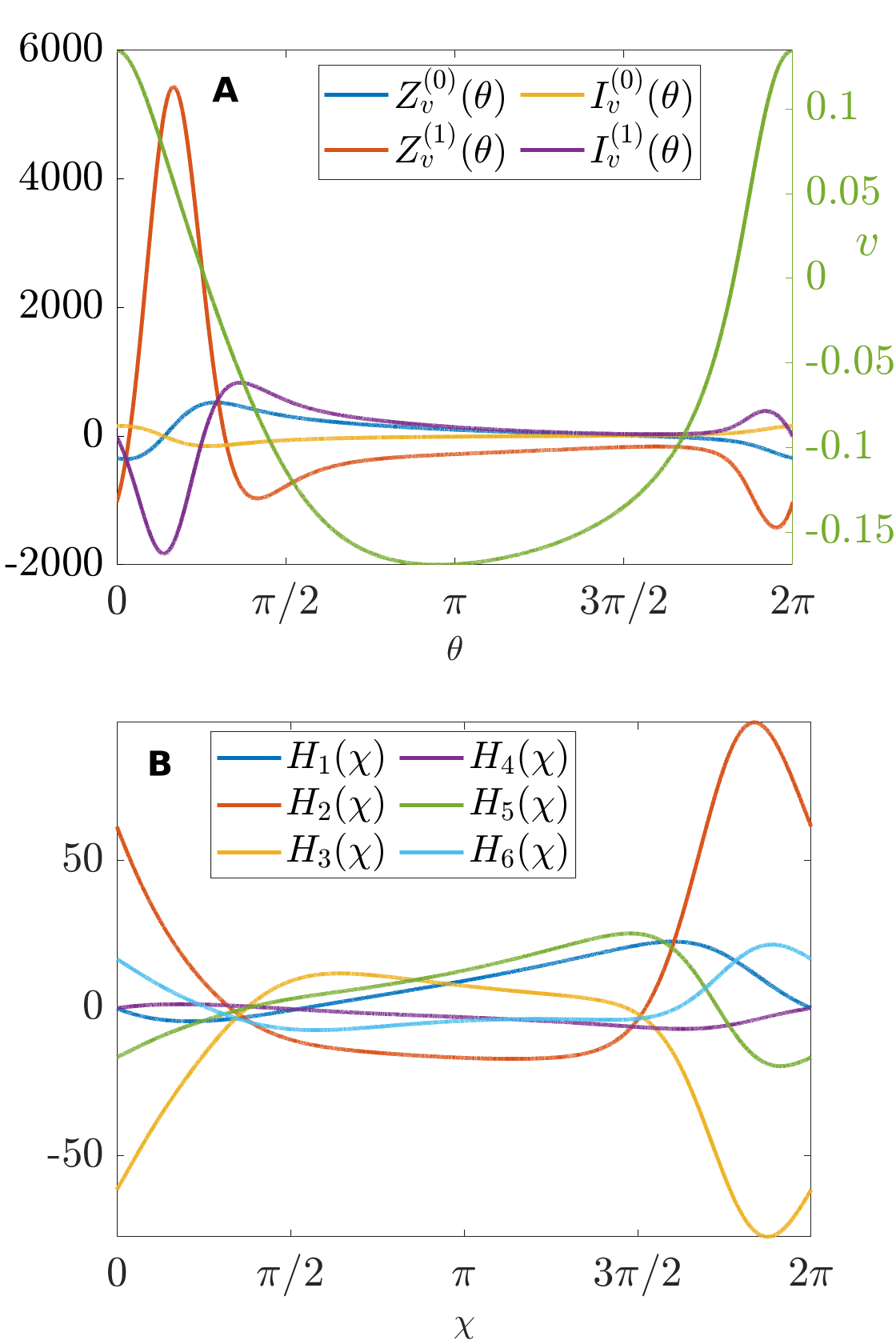}
\caption{(A) The first (voltage) components of response  vectors $Z^{(0)}$, $I^{(0)}$, $Z^{(1)}$ and $I^{(1)}$ together with the voltage component $v(\theta)$ of the limit cycle,  and (B) the interaction functions $H_1, \ldots, H_6$ given by \eqref{eq:Hk} for the Morris-Lecar model \eqref{eqs:ML} with parameter values as in Appendix \ref{sec:AppendixE} and linear coupling.
 }
\label{fig:MLfuncs}
\end{center}
\end{figure}

The Floquet eigenfunction $g^{(1)}$, and the iPRC and iIRC, $Z^{(0)}$ and $I^{(0)}$ and their first order correction terms $Z^{(1)}$ and $I^{(1)}$ are numerically determined as appropriately normalised solutions of \eqref{eqs:adjoints} and the first (voltage) component of the response vectors are depicted together with the limit cycle in Figure~\ref{fig:MLfuncs}A. 

We consider networks diffusively coupled through the voltage variable $v$ so that $\vec{x}_i = (v_i, w_i)^T$ and $\vec{G}(\vec{x}_i, \vec{x}_j) = (v_j-v_i, 0)^T$ and $w_{ij}=1/N$ in \eqref{eq:fullnetwork}. Figure~\ref{fig:MLfuncs}B shows the resulting six interaction functions $H_1, \ldots, H_6$ given by \eqref{eq:Hk}. In section \ref{sec:ML2} we compare the dynamics of the phase-isostable approximation with that of the full model for two coupled nodes and consider a larger network of 200 nodes in section \ref{sec:ML200}. 

\subsubsection{Two coupled nodes}\label{sec:ML2}

\begin{figure*}
\begin{center}
\includegraphics[width = \textwidth]{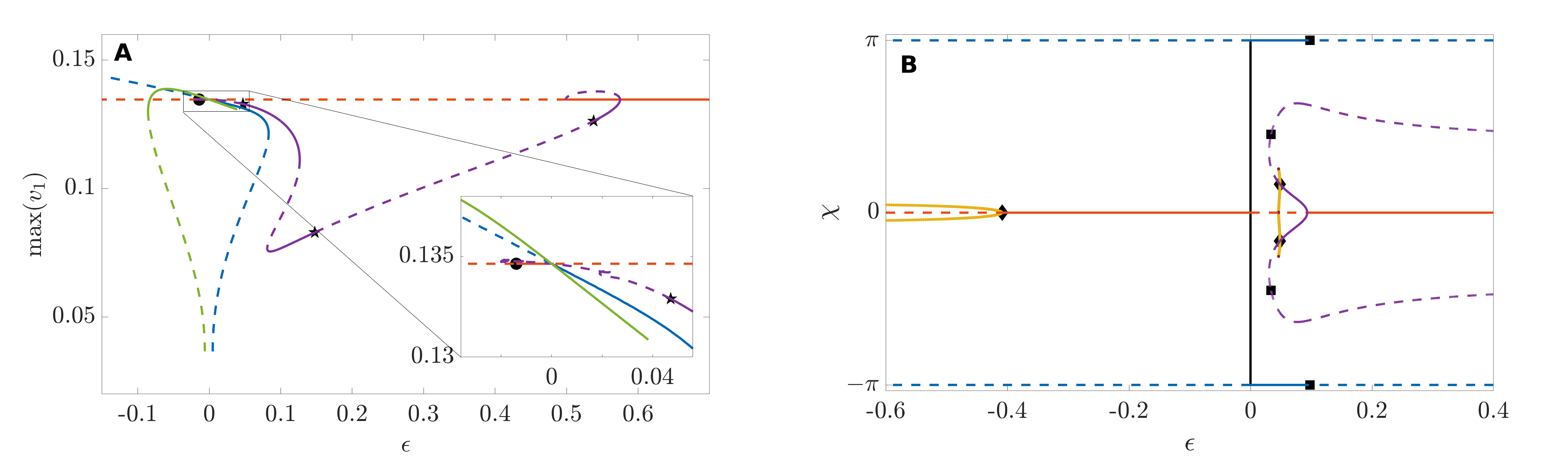}
\caption{Bifurcation diagrams for two linearly coupled Morris-Lecar neurons with interaction strength $\epsilon$ as the bifurcation parameter. (A) shows the full model equations \eqref{eqs:ML} with the maximum of the voltage variable for one of the two nodes on the vertical axis. (B) shows the approximations using the phase-isostable network equation framework \eqref{eqs:piaveraged} with the phase difference $\chi= \theta_2- \theta_1$ on the vertical axis. Stable (unstable) periodic solutions are shown by solid (dashed) lines with synchrony in light red, antisynchrony in blue, phase locked states with $\chi\neq 0, \pi$ in purple and other periodic states in green ((A) only). In (B), stable (unstable) quasiperiodic solutions are indicated in yellow (dark red). Black circles/stars/diamonds indicate period doubling/torus/Hopf bifurcations respectively (noting that periodic solutions are given by equilibrium points of the phase-isostable approximation when working with the phase difference $\chi$.) Black squares show limit points where the modulus of an isostable coordinate asymptotes to $\pm\infty$. The phase-isostable approximation captures many of the qualitative features of the full model dynamics as discussed in the text. Parameter values for the node dyanamics for both diagrams are as given in Appendix \ref{sec:AppendixE}. The bifurcation diagrams are produced using the implementation of AUTO in XPPAUT \cite{XPPAUT}. }
\label{fig:MLbif}
\end{center}
\end{figure*}

For a network of two nodes, bifurcation diagrams using the interaction strength $\epsilon$ as the bifurcation parameter are computed for both the full model (Figure~\ref{fig:MLbif}A) and the phase-isostable network equations \eqref{eqs:piaveraged} (Figure~\ref{fig:MLbif}B). For the full model the bifurcation diagram, Figure~\ref{fig:MLbif}A, shows the maximum of the voltage variable for one of the two nodes on the vertical axis for relatively small values of the coupling strength $\epsilon$. In the corresponding diagram for the phase-isostable approximation (Figure~\ref{fig:MLbif}B) the vertical axis shows $\chi = \theta_2-\theta_1$, the phase difference between the two nodes, so that periodic orbits of the full model correspond to equilibrium points of the phase-isostable reduced system. Since we choose to work with phase-isostable network equations \eqref{eqs:piaveraged} that are linear in the isostable coordinates, we find unique branches of synchronous and antisynchronous solutions with the isostable coordinate given by \eqref{eq:synciso} and \eqref{eq:finiteNsplay} respectively in Figure~\ref{fig:MLbif}B. These branches also exist for all values of $\epsilon$, except where the value of the isostable coordinate asymptotes to positive or negative infinity. There may be more than one phase-locked state with $\chi \neq 0, \pi$ for a given value of $\epsilon$ depending on the interaction functions $H_1, \ldots, H_6$. The existence, stability of solution branches and locations of bifurcations in Figure~\ref{fig:MLbif}B all agree with the explicit calculations in Section \ref{sec:phaselocked}.  Taking Figure~\ref{fig:MLbif}A as ground truth, we here comment on the extent to which the phase-isostable framework approximates the dynamics observed in the full system as the strength $\epsilon$ of the linear coupling increases. 

Considering first the synchronous state (in light red) we see that in the full model this solution is unstable in the weak coupling regime for positive $\epsilon$, but restabilises at $\epsilon= 0.499$ where it meets the branch of phase-locked solutions (shown in purple). Synchrony is also stable for a small interval of negative $\epsilon$, before undergoing a period doubling bifurcation at $\epsilon =-0.015$. The phase-isostable approximation, as shown in Figure~\ref{fig:MLbif}B, correctly predicts the stability type of synchrony in the neighbourhood of $\epsilon=0$ (as does a first-order phase approximation). The phase-isostable approximation is additionally able to capture bifurcations of synchrony for both positive and negative $\epsilon$, showing synchrony gaining stability where it meets a phase-locked state at $\epsilon =0.0934$ and losing stability at a Hopf bifurcation at $\epsilon =-0.407$. We note that while this reproduces the qualitative behaviour of the synchronous branch, the phase-isostable approximation underestimates the precise value of $\epsilon$ at which the bifurcations occur, just as for the MF-CGLE in section \ref{sec:MFCGL}. 

In the full system antisynchrony (shown in blue in Figure~\ref{fig:MLbif}) is unstable for all negative values of $\epsilon$, but is stable for $0<\epsilon<0.0831$ where the orbit lies inside, but close to the synchronous orbit $\gamma$. For $0.00491 < \epsilon< 0.0831$ the stable large amplitude antisynchronous solution co-exists with an unstable smaller amplitude antisynchronous solution. The qualitative stability of the branch of antisynchronous solutions of the full model with orbit closest to the uncoupled node orbit is also captured by the phase-isostable reduction, including the fact that it exists for $\epsilon<\epsilon_\infty$ where the phase isostable framework estimates that $\epsilon_\infty = 0.0961$ (compared to $0.0831$ in the full system). This point is labelled as a limit point (black square) in Figure \ref{fig:MLbif}. As $\epsilon \to \epsilon_\infty$ from below the stable solution branch has $\Psi\rightarrow -\infty$ corresponding to a shrinking amplitude orbit in the $(v,w)$ coordinates approaching the unstable fixed point inside the uncoupled node limit cycle. We also find a branch of unstable antisynchronous solutions for $\epsilon>\epsilon_\infty$ with $\Psi$ positive. As $\epsilon \to \epsilon_\infty$ from above this branch has $\Psi\rightarrow\infty$ (corresponding to a growing amplitude orbit approaching the outer edge of the basin of attraction of the uncoupled node limit cycle, see Figure \ref{fig:MLphaseplane}). There is no bifurcation at $\epsilon_\infty$; the change in stability here is coincidental.  We have not found a corresponding antisynchronous state with orbit outside of $\gamma$ in the full model since it may always be unstable. The phase-isostable reduction at the given order does not find the unstable small amplitude antisynchronous oscillations since, as previously noted, it gives a unique value of ${\Psi}$ for each $\epsilon$. It may be that if the equations \eqref{eqs:piaveraged} were taken to higher order in $\epsilon$, then this other branch could be revealed. We note that the phase-isostable framework predicts a small region ($0.0934<\epsilon<\epsilon_\infty$) where both synchrony and antisynchrony are stable. This does not occur in the full system. 

\begin{figure}
\begin{center}
\includegraphics[width = 0.45\textwidth]{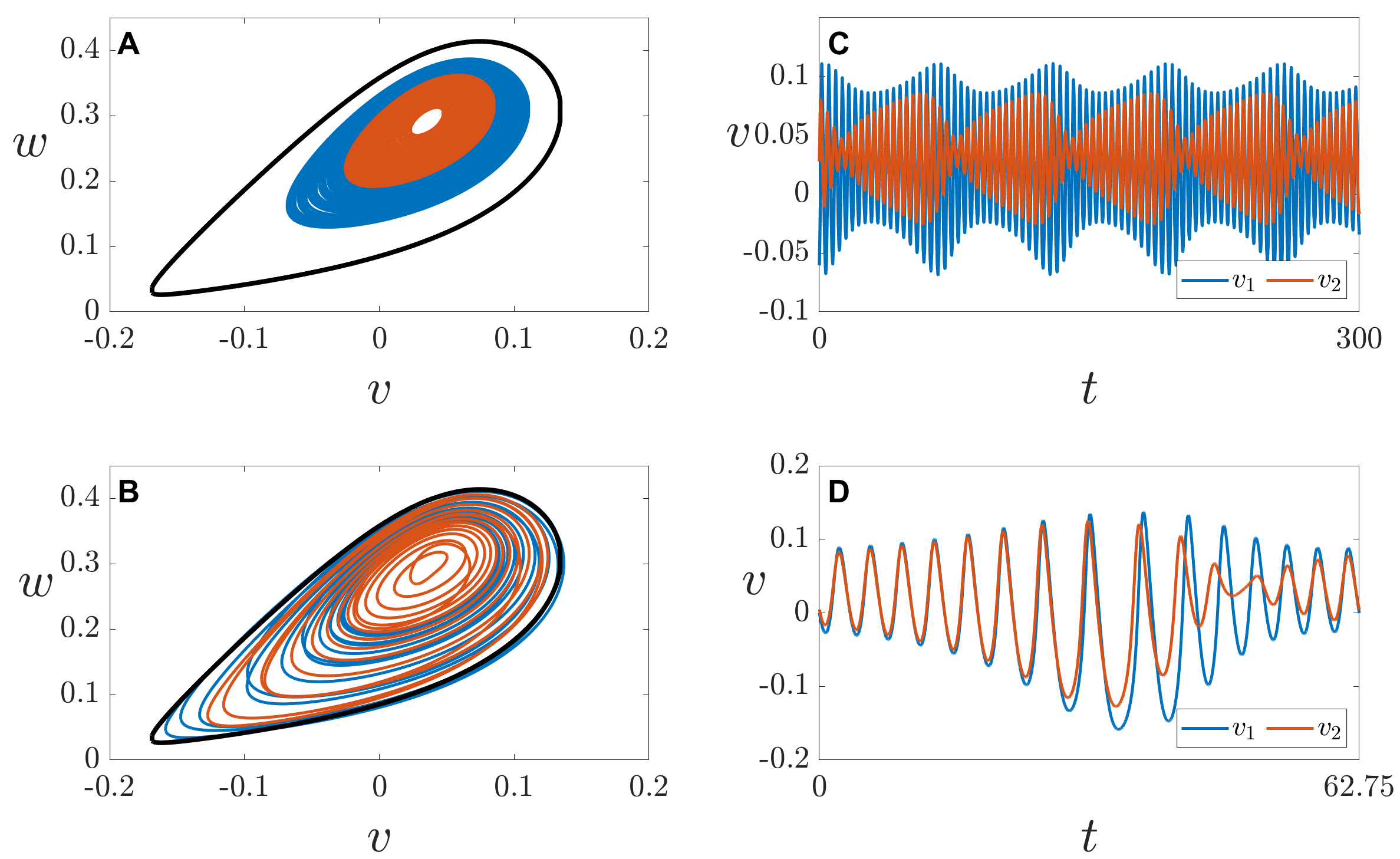}
\caption{Direct numerical simulations of two linearly coupled Morris Lecar neurons for values of $\epsilon$ where all solutions indicated in Figure~\ref{fig:MLbif}A are unstable. (A,C) $\epsilon=0.25$, shows quasiperiodic behaviour for both nodes, as is typical for $\epsilon \in (0.148, 0.375)\cup (0.403, 0.499)$ (B,D) $\epsilon =0.39$, shows periodic behaviour for both nodes, as is typical for $\epsilon \in (0.375, 0.403)$. Panels (A,B) show network activity in the $(v,w)$ plane (synchronous orbit $\gamma$ is indicated in black), while the corresponding time series for $v_1, v_2$ are shown in (C,D). In $D$ a single period is shown.  Node parameters are as in Appendix \ref{sec:AppendixE}.}
\label{fig:MLquasi}
\end{center}
\end{figure}

In the full coupled system there are also periodic solutions for which the nodes do not share a common orbit and therefore for each of these there is a symmetric solution under $\vec{x}_1 \leftrightarrow \vec{x}_2$. For clarity, Figure~\ref{fig:MLbif}A shows only the branch for which $v_1$ has the largest maximum value. A phase-locked solution with $\vec{x}_1$ performing larger and $\vec{x}_2$ smaller amplitude oscillations inside the synchronous orbit exists for $-0.0218 <\epsilon< 0.575$ (shown in purple in Figure~\ref{fig:MLbif}). Its stability varies along the branch as it undergoes a series of saddle node and torus bifurcations. Many of these occur within the region near $\epsilon=0$, but most obvious in Figure~\ref{fig:MLbif}A are the torus bifurcations at $\epsilon = 0.0472$, $\epsilon=0.148$ and $\epsilon=0.538$ marked by black stars. There is also an isola of  periodic solutions for $0.375<\epsilon< 0.403$ an example of which is depicted in Figure~\ref{fig:MLquasi}B,D. This solution is not shown in Figure~\ref{fig:MLbif}A to avoid confusion with the synchronous branch since they share similar values of $\max(v_1)$. For values of $\epsilon \in (0.148, 0.375)\cup (0.403, 0.499)$ numerical simulations indicate that the stable solutions are quasiperiodic with oscillations fluctuating around the phase-locked state, as indicated in Figure~\ref{fig:MLquasi}A,C. In the phase-isostable approximation, branches of stable phase-locked states with $\chi \neq 0, \pi$ bifurcate from synchrony at $\epsilon = 0.0934$ losing stability at a Hopf bifurcation at $\epsilon=0.0484$. This gives rise to stable solutions where $\chi, \Psi_1, \Psi_2$ are all periodic corresponding to behaviour observed in numerical simulations of the full model (Figure~\ref{fig:MLquasi}), however, here these solutions exist over a much smaller interval of $\epsilon$ values. We also note that the stable phase-locked state has isostable coordinates of opposite sign for the two nodes corresponding to one node orbiting outside of $\gamma$ and the other inside. As the limit point at $\epsilon = 0.0339$ is approached, one of the isostable coordinates asymptotes to positive infinity while the other asymptotes to negative infinity. Beyond the limit point there is another (unstable) branch of phase-locked states with the signs of the isostable coordinates reversed. This branch of solutions persists for larger values of $\epsilon$, but is unstable and always has one node with a large isostable coordinate.

The branch of periodic solutions shown in green in Figure~\ref{fig:MLbif}A is created at a homoclinic bifurcation at $\epsilon=0.0381$ and meets the unstable fixed point at $\epsilon=-0.0618$. This solution has $\vec{x}_1$ oscillating near the synchronous orbit $\gamma$, while $\vec{x}_2$ performs very small amplitude oscillations near the fixed point outside of the basin of attraction of $\gamma$. Since phase-isostable coordinates are valid only within the basin of attraction of $\gamma$, the phase-isostable network equations are not able to capture this oscillatory solution.

We conclude that the phase-isostable approximation is able to capture the qualitative stability of synchrony and antisynchrony away from $\epsilon=0$. Furthermore, it is able to reveal phase-locked states and quasiperiodic solutions, capturing far more of the full dynamics than the first-order phase reduction of the model. Since approximation using phase-isostable coordinates can describe qualitative dynamics which first-order phase reduction cannot for a two node network, we now use the approximation to investigate the dynamics of larger networks of Morris-Lecar models in section \ref{sec:ML200}. We note that, while the phase-isostable framework can be instructive concerning network behaviours for relatively weak coupling, since the expansions rely on the assumption that $\psi = O(\epsilon)$ predictions for larger values of $\epsilon$ must be interpreted with caution.

\subsubsection{Networks of many Morris-Lecar neurons}\label{sec:ML200}

\begin{figure}
\begin{center}
\includegraphics[width = 0.45\textwidth]{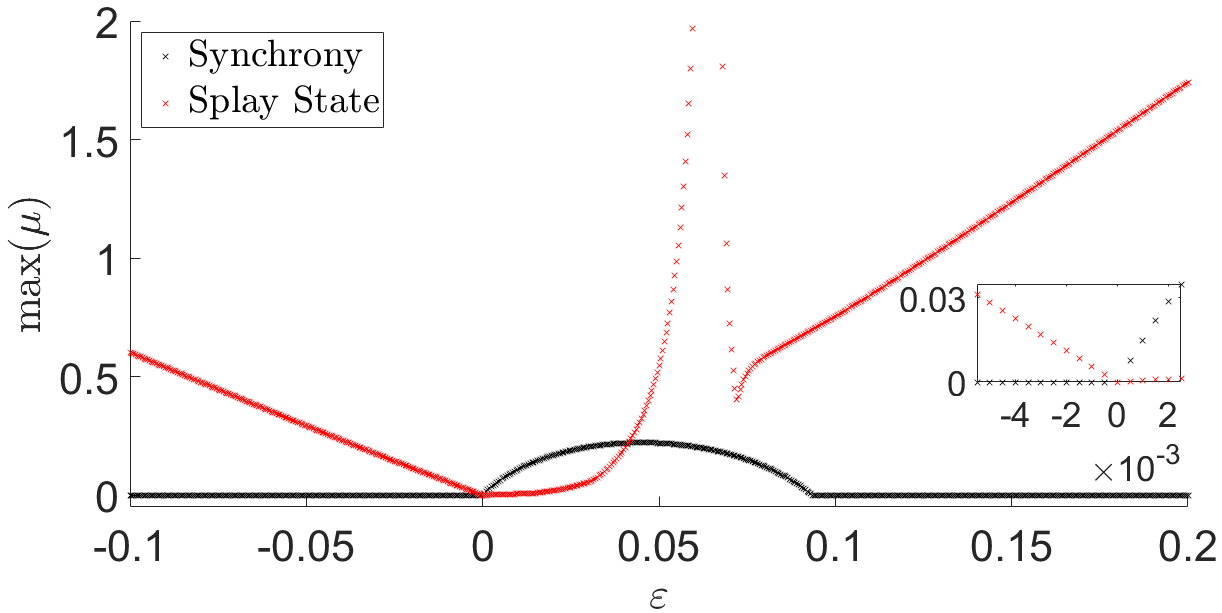}
\caption{The maximal real part of eigenvalues of phase-locked states in a 200 node, globally-coupled, phase-isostable Morris-Lecar network plotted against the coupling strength $\epsilon$. The synchronous state is shown in black and the splay state is shown in red.}
\label{fig:ML200NodeEigenval}
\end{center}
\end{figure}

We next consider a network of 200 diffusively coupled Morris-Lecar neurons. We again use the parameters values in Appendix \ref{sec:AppendixE} so that the response and interaction functions remain as in Figure \ref{fig:MLfuncs}. Using \eqref{eq:stabsyncpi} and \eqref{eq:splayeigs} we calculate the stability of the synchronous and splay states respectively over a range of small positive and negative coupling strength $\epsilon$. The results are shown in Figure~\ref{fig:ML200NodeEigenval}. The maximal real part of the eigenvalues for synchrony is positive (and therefore the state is unstable) when $0<\epsilon<0.0934$ and the splay state is unstable for all values of $\epsilon\neq 0$. The discontinuity in the maximum real part of the eigenvalues for the splay state at $\epsilon=0.0664$ corresponds to the discontinuity in the isostable coordinate of the splay state when the denominator in \eqref{eq:finiteNsplay} becomes zero. The splay state orbit lies inside the node orbit ($\Psi<0$) when $\epsilon<0.0664$ and for $\epsilon>0.0664$ the splay state has $\Psi>0$. 

\begin{figure}
\begin{center}
\includegraphics[width = 0.4\textwidth]{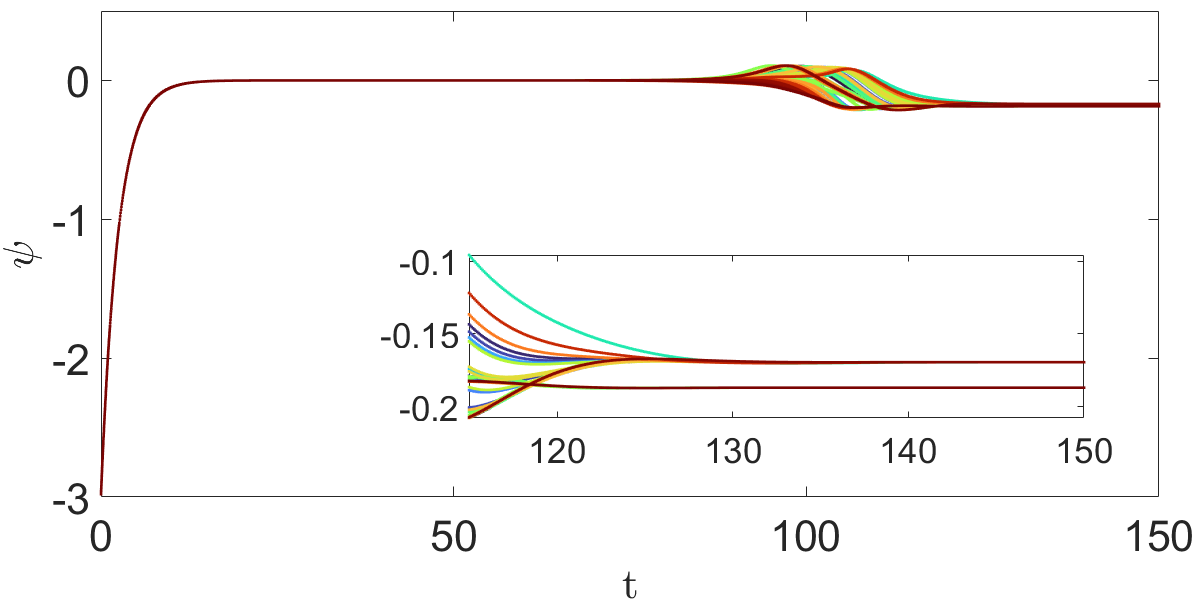}
\caption{The isostable coordinates of the 200 nodes in the phase-isostable approximation of the diffusively coupled network of Morris-Lecar neurons evolving with time. The initially near synchronous group moves towards the synchronous isostable where $\psi_i=0$ for all $i$ before desynchronising at $t\approx 80$ and settling into a stable two cluster state from $t \approx 130$. In this cluster state, there are 28 nodes in the cluster with $\psi  =-0.1694$ and 172 nodes in the cluster with $\psi=-0.1868$. Interaction functions are as in Figure~\ref{fig:MLfuncs} and $\epsilon =0.065$. }
\label{fig:200PhaseIso}
\end{center}
\end{figure}

We take $\epsilon = 0.065$ which lies in the range where both synchrony and the splay state in the phase-isostable approximated framework are unstable and numerically simulate the phase-isostable network equations to investigate the stable behaviour predicted in this region. We initialise the system with the 200 nodes in a tightly packed group with phase coordinates $\theta_i \in [0.283725, 0.283735]$ and isostable coordinates $\psi_i \in [2.9794, 2.9798]$. Figure~\ref{fig:200PhaseIso} shows the evolution in time of the isostable coordinates. We see the group of nodes initially move towards the isostable $\psi=0$ and remains near the synchronous network state for some time. The group then rapidly desynchronises before settling into a stable 2-cluster state. The analysis in section \ref{sec:twoclusters} confirm the existence of the observed cluster state with $N_A=28$ and $N_B=172$ and also that it has phase difference $\chi =2.1407$ between the clusters which have orbits coinciding with the $\Psi_A=-0.1694$ and $\Psi_B=-0.1868$ isostables, also in agreement with the simulation. The stability analysis further confirms that this state is stable for $\epsilon=0.065$.  

\begin{figure}
\begin{center}
\includegraphics[width = 0.4\textwidth]{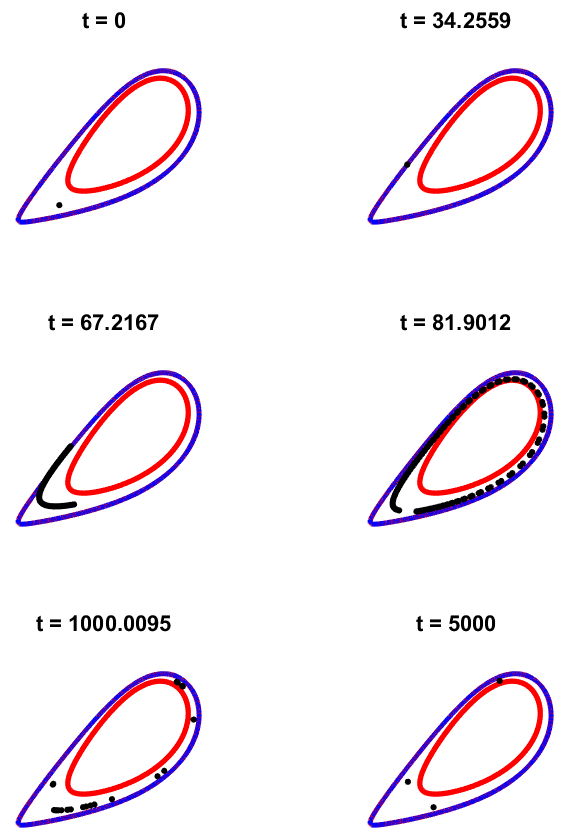}
\caption{A series of snapshots of a simulation of the full network of 200 diffusively coupled Morris-Lecar neurons showing positions of each node (black dots), the uncoupled limit cycle at $\psi=0$ (blue) and the splay state orbit at $\psi=-19.3$ (red). Parameter values are as in Appendix \ref{sec:AppendixE} with $\epsilon =0.065$. }
\label{fig:200NodeSnapshots}
\end{center}
\end{figure}

For the same value of the coupling strength $\epsilon = 0.065$ we also simulate the full network equations from an initial state with $(v_i, w_i)$ in a group close to $(\overline{v}, \overline{w}) \approx $ (-0.1,0.07) which has phase and isostable coordinates $(\theta, \psi) = (0.28373,2.9796)$. Snapshots from the simulation are shown in Figure~\ref{fig:200NodeSnapshots}, which also indicates the isostable $\psi=0$ which coincides with the synchronous orbit and the isostable corresponding to the splay state for these parameter values ($\psi =-19.3$) in the $(v,w)$ space. Initially the nodes are close to synchronous and are first drawn towards and orbit near to $\psi=0$. However, since this state is unstable the cluster soon breaks up and the nodes spread out and oscillate nearer to the splay state isostable. This state is also unstable in the full network and therefore the system oscillates between the splay and synchronous states for a time before settling at a stable 3-cluster state. The simulations of the phase-isostable approximation and the full model for equivalent initial conditions show the same progression from near synchrony, subsequent desynchronisation eventually settling on a cluster state. 

For stronger diffusive coupling ($\epsilon =0.4$) with the node parameters still as in Appendix \ref{sec:AppendixE}, Han et al. \cite{Han1995} observed that in the 200 node network initialised at a near synchronous state, the nodes would at first behave as we have observed at lower coupling strength, moving to the synchronous orbit before desynchronising. However at this larger coupling strength it was observed that the desynchronised nodes spiral in towards the unstable node inside the limit cycle, before again moving out towards the synchronous orbit. This behaviour repeats resulting in mean field voltage traces showing large amplitude fluctuations. Unfortunately this type of behaviour cannot be captured using the phase-isostable network equations at the order of \eqref{eqs:piaveraged} as the regime where $\epsilon=0.4$ appears to be beyond the scope of its predictive power. Nonetheless, we have been able to demonstrate that at smaller coupling strengths the phase-isostable framework is able to capture the formation of stable cluster states for parameter values where both the synchronous and splay states are unstable in agreement with simulations of the full network equations. 

\subsection{Remarks}
\begin{figure}
\begin{center}
\includegraphics[width = 0.5\textwidth]{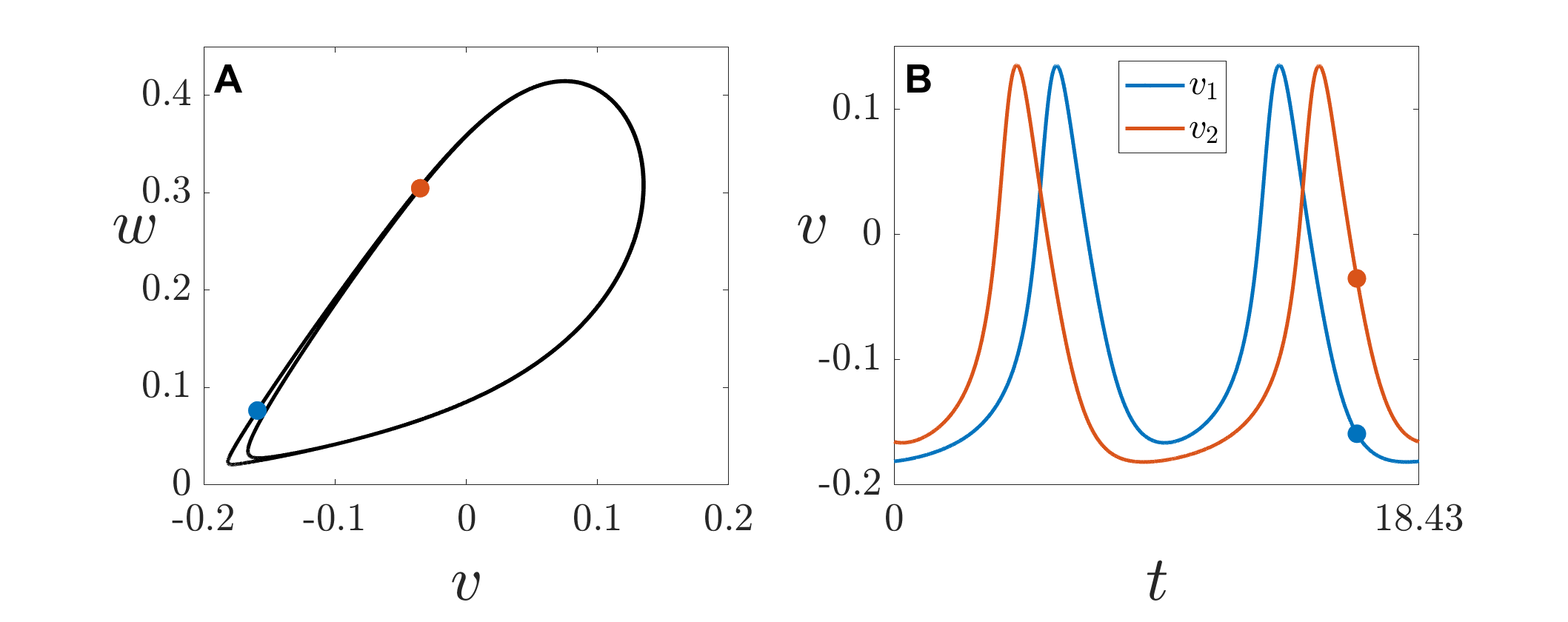}
\caption{(A) A stable periodic phase-locked state in the network of two coupled Morris-Lecar neurons occurring at $\epsilon = -0.0129$ on the solution branch bifurcating from the periodic doubling bifurcation indicated at $\epsilon=-0.015$ in Figure \ref{fig:MLbif} A. Since this periodic orbit projected into the two-dimensional node phase space has a self intersection, this orbit cannot be described by the linear phase-isostable network equations whose phase-locked solutions all coincide with isostables which do not have self intersections. (B) shows the corresponding time series for $v_1$, $v_2$ over a single period.}
\label{fig:MLtrajpd}
\end{center}
\end{figure}

In section \ref{sec:phaselocked} we observed that in the phase-isostable network equations truncated at linear order in the isostable coordinates (or equivalently $O(\epsilon^2)$) \eqref{eqs:piaveraged}, the orbits for individual nodes all coincide with isostables since they have constant isostable coordinate in the most slowly decaying direction (which is the only direction we retain here). When the node dynamics are two-dimensional the phase-isostable coordinates are a transformation from the original coordinates of the node model. That is, for each $\vec{x} \in \mathcal{B}(\gamma)$ there is a unique corresponding $(\theta, \psi)$. Therefore, since trajectories of the node dynamics cannot intersect, isostables also have no self intersections. However, for periodic orbits of the full network dynamics the projection of the dynamics of each node onto the node phase space can have shared orbits which intersect. For example, Figure \ref{fig:MLbif}A indicates a period doubling bifurcation of synchrony at $\epsilon= -0.015$. The period doubled branch (not shown in Figure \ref{fig:MLbif}A) at $\epsilon =-0.0129$ is a periodic two-cluster phase-locked state which is bistable with synchrony. The shared node orbit is shown in Figure \ref{fig:MLtrajpd}A and we observe that this has a self-intersection. The linear order truncation of the phase-isostable network equations used here cannot capture such solutions, but it is possible that a higher order expansion would be capable of recovering such orbits. 
Higher order truncation in the network equations \eqref{eqs:piaveraged} would also likely improve the quantitative accuracy of the predictions of the phase-isostable dynamics with the dynamics of the full system. 

\section{Conclusions}\label{sec:conclusions}

There are two emerging frameworks providing improved understanding of the behaviour of coupled oscillator networks through the use of isostable coordinates to capture dynamics off-cycle. While many authors are choosing to use isostable dynamics to obtain higher-order phase reductions \cite{Mau2023, Park2021, Park2023, Wilson2019c} we instead focus on the approach which retains the isostable dynamics in the slowest decaying direction. We have obtained existence and stability criteria for phase-locked states in networks of identical nodes and used these to demonstrate that for the MF-CGLE the phase-isostable network equations have far superior capabilities at capturing the qualitative shape of the stability boundaries of the synchronous and splay states in networks of arbitrary size when compared with phase-reduction at up to cubic order in the interaction strength. For a network of size $N=3$ the phase-isostable framework can also correctly identify the stable solution types beyond the loss of stability of synchrony and the splay state. We then use the framework to study networks of 2 and 200 diffusively coupled Morris-Lecar neurons. For the smaller network, comparison of bifurcation diagrams for both the full and phase-isostable approximations of the dynamics show that the phase-isostable framework can capture qualitative changes in stability of synchrony and the antisynchrnous state as well as the existence of other phase-locked states and quasiperiodic solutions which cannot be described by first-order phase reduction. However, it is not able to accurately predict the values of the coupling strength at which bifurcations occur, tending to underestimate these values. For a network of 200 nodes, for values of the coupling strength where both synchrony and the splay state are unstable in the full system and the phase-isostable approximation, the phase isostable network equations predict stable cluster states which are also seen in numerical simulations of the full equations in this regime. The phase-isostable network approximation (or reduction in the case of node dynamics of dimension three or higher) appears to be a useful tool to indicate qualitative network behaviour beyond that which can be revealed using first-order phase reduction whilst also keeping the computational complexity low since only six interaction functions need to be computed, requiring only four response functions.

For expansions of comparable order and computational effort, we have observed that retaining the notion of distance from cycle through the isostable dynamics outperforms using the isostable dynamics to refine phase-reduced dynamics in locating bifurcations of phase-locked solutions. Furthermore it allows for the analysis of off cycle dynamics through the tracking of the isostable coordinate.  
In section \ref{sec:MFCGL} we found that the phase-isostable network transformation best captures the qualitative bifurcation structure of the MF-CGLE. Since it is the normal form for globally diffusively coupled Hopf normal form oscillators, it is natural to suppose that this framework is the most appropriate to use to study linearly coupled networks of oscillatory nodes near Hopf bifurcation. We do however note that recent work of Bick et al. \cite{Bick2023} generalises the techniques of \cite{Leon2019} to consider Stuart-Landau like systems with phase-dependent amplitude emulating the deformed limit cycles expected away from Hopf bifurcation. They show that a second-order phase reduction is able to accurately predict the stability properties of the synchronized and splay orbits when all terms of up to second order in the coupling strength are included. It would be interesting to compare with corresponding results using a phase-isostable approach. 

Although we have chosen to use the orbit of an uncoupled oscillatory node as the reference for setting up the phase-amplitude coordinate system, other choices are possible.  This has been recognised by Wilson as a way to to handle perturbations that take one far from cycle with the development of an adaptive reduction strategy to construct a more useful reference limit cycle (so that distances to this reference cycle remain small) \cite{Wilson2020a}.
Alternatively, one could more simply take a phase-isostable reduction about a stable periodic network state (e.g. synchrony or another phase-locked state) \cite{Wilson2019b}. This approach has been used to investigate control strategies for desynchronising populations of neural oscillators under periodic stimulation \cite{Wilson2020c, Wilson2022}.

In the present work we have investigated finite size networks with linear (state-dependent) coupling and two-dimensional node dynamics. There are obvious extensions to this work, incorporating other forms of coupling such as event driven interactions, along the lines described for synaptic interactions e.g. in  \cite{Nicks2018}, and also to investigate networks with higher dimensional node dynamics where retaining a single isostable coordinate dynamics in the most slowly decaying direction represents a dimension reduction for the system.  The extension to treat dense graphs could naturally build on the work in \cite{Medvedev2019,Ihara2023} for Kuramoto networks using the notion of a graphon. Moreover, the analysis of continuum models with non-local interactions using phase-amplitude coordinates seems feasible by generalising the approach in \cite{Wilson2019c} developed for reaction--diffusion equations.  Finally, let us mention the challenge of dealing with delays.  At the node level these can induce oscillations, and at the network level can strongly influence patterns of phase-locked states and their bifurcations.  A method for constructing the infinitesimal phase response for delay induced oscillations has already been developed in \cite{Kotani2012,Novicenko2012}, and it would be interesting to use this functional setting to develop the corresponding amplitude response, as well as to incorporate delayed interactions within a phase-amplitude network setting. 
All of the above are topics of ongoing work and will be reported upon elsewhere.


\appendix

\section{Higher-order terms in PRC and IRC expansions} \label{sec:AppendixA}

When the $n$ dimensional system \eqref{eq:uncoupled} has been reduced to the dynamics of the phase $\theta$ and a single isostable coordinate $\psi$ as in \eqref{eq:phaseiso}--\eqref{eq:phaseiso2}, the expansion to any order of the solution about the periodic orbit, the PRC and the IRC are given by \eqref{eq:expg}--\eqref{eq:expI}. Following \cite{Wilson2020b}, and simplifying for the case of a single isostable coordinate we here derive the adjoint equations \eqref{eq:eigODE}--\eqref{eq:IRCODE} which must be satisfied by the vector functions $g^{(k)}$, $Z^{(k)}$ and $I^{(k)}$.

\subsection{Asymptotic expansion of eigenfunctions}

As in \eqref{eq:expg} the deviation of a trajectory from the point with the same phase on the limit cycle may be expressed as
\begin{align}\label{eq:expdeltax}
\Delta \vec{x}(\theta, \psi) = \sum_{k=1}^\infty \psi^kg^{(k)}(\theta).
\end{align} Differentiating with respect to $t$ and using the fact that $\dot{\psi}= \kappa \psi$ we have
\begin{align}
\label{eq:diffdeltax}
\FD{\Delta \vec{x}}{t} = \sum_{k=1}^\infty \psi^k \left( \FD{g^{(k)}(\theta)}{t}+ k\kappa g^{(k)}(\theta)\right).
\end{align}
With the notation $\vec{F}(\vec{x}) = [F_1(\vec{x})  \ldots  F_n(\vec{x})]^T$, matrices of partial derivatives evaluated at $\vec{x}^\gamma(\theta)$ on cycle can be defined recursively by
\begin{align}\label{eq:Frecurrive}
F_q^{(k)}(\theta) = \PD{\mbox{vec}(F_q^{(k-1)})}{\vec{x}^T} \in \mathbb{R}^{(k-1)n \times n},
\end{align} where $F_q^{(0)}(\theta) = F_q(\vec{x}^\gamma (\theta))$. Here $\mbox{vec}(\cdot)$ is an operator which stacks each column of a matrix resulting in a single column vector \cite{Wilson2020b, Magnus2007}. The Taylor expansion of \eqref{eq:uncoupled} about $\vec{x}^\gamma(\theta)$ then gives
\begin{align}\label{eq:diffdeltax2}
\FD{\Delta \vec{x}}{t} = J\Delta\vec{x} + \begin{bmatrix}
\sum_{i=2}^{\infty} \frac{1}{i!}[\Delta \vec{x}^{T}]^{\otimes i}\mbox{vec}(F_{1}^{(i)}(\theta))\\
\vdots \\
\sum_{i=2}^{\infty} \frac{1}{i!}[\Delta \vec{x}^{T}]^{\otimes{i}}\mbox{vec}(F_{n}^{(i)}(\theta))\\
\end{bmatrix},
\end{align} where $[\Delta \vec{x}^{T}]^{\otimes i} = \Delta \vec{x}^{T} \otimes \cdots \otimes \Delta \vec{x}^{T}$ with $\Delta \vec{x}^{T}$ appearing $i$ times and $\otimes$ is the Kronecker product. Upon substituting \eqref{eq:expdeltax} into \eqref{eq:diffdeltax2}, equating with \eqref{eq:diffdeltax} and matching terms in $\psi^k$ we obtain a first order differential equation for $g^{(k)}$ in terms of $g^{(j)}$ where $1<j<k$ \cite{Wilson2020b}. We can obtain an explicit expression for these differential equations by noting that
\begin{align}\label{eq:kroneckerDeltax}
[\Delta \vec{x}^{T}]^{\otimes i} = \left[\sum_{k=1}^{\infty} \psi^{k} g^{(k)T} \right]^{\otimes i} = \left[ \sum_{j=i}^{\infty} \psi^{j} \sum_{\mathcal{A}^{i,j}} \bigotimes_{p = 1}^{i} g^{(a_{p})}\right]^T
\end{align} where $\mathcal{A}^{i,j}$ is the set of all ordered lists $(a_1, a_2 \ldots a_i)$ such that $a_i \in \mathbb{N}$ and $\sum_{p=1}^ia_p =j$ (for example $\mathcal{A}_{2,3} = \{ (1,2), (2,1)\}$ ) and for a list $(a_1, a_2 \ldots a_i) \in \mathcal{A}_{i,j}$,
\begin{align}
\bigotimes_{p = 1}^{i} g^{(a_{p})} = g^{(a_1)} \otimes \cdots \otimes g^{(a_i)}.
\end{align}  We then find that matching powers of $\psi^k$ results in the ODE \eqref{eq:eigODE} for $g^{(k)}$ where $\alpha^{(k)} = \left[ \alpha_{1}^{(k)}, \ldots  \alpha_{n}^{(k)}\right]^T$ with
\begin{align}\label{eq:alpha}
\alpha^{(k)}_i &= \sum_{l=2}^k \frac{1}{l!} \left[ \sum_{\mathcal{A}^{l,k}}\bigotimes_{p = 1}^{l} g^{(a_{p})} \right]^T \mbox{vec}(F_i^{l}(\theta)), \quad i=1, \ldots, n.
\end{align} Thus $\alpha^{(k)}$ is composed of only the lower order terms $g^{(j)}$, $1\leq j <k$ and therefore the expansion \eqref{eq:expg} can be obtained by solving the equations from the lowest order to successively higher order. It is noted in \cite{Wilson2020b} that the equation for $g^{(1)}$ requires normalisation since it has a zero Floquet exponent. We choose the normalisation $|g^{(1)}(0)|=1$ in contrast to \cite{Wilson2020b} who take $I^{(0)}(0)\cdot g^{(1)}(0) =1$ as both the normalisation for $g^{(1)}$ and $I^{(0)}$ leading to non-uniqueness of the expansions.

\subsection{Asymptotic expansion of the PRC and IRC}

Again following Wilson \cite{Wilson2020b}, we now consider the expansions of $\mathcal{Z}(\theta, \psi)$ and $\mathcal{I}(\theta, \psi)$ of \eqref{eq:expZ} and \eqref{eq:expI} respectively. Using standard arguments, Wilson \cite{Wilson2020b} shows that $\mathcal{Z}$ satisfies the adjoint equation
\begin{align}\label{eq:Zadjoint}
\FD{\mathcal{Z}}{t} = -\left.\PD{\vec{F}^T}{\vec{x}}\right|_{\vec{x}} \mathcal{Z}
\end{align} which differs from \eqref{eq:PRCadjoint} since it describes the phase response to perturbations of a trajectory which is already away from the limit cycle. Then we may make the asymptotic expansion
\begin{align}
&\left.\PD{\vec{F}^T}{\vec{x}}\right|_{\vec{x}^\gamma + \Delta \vec{x}} = J^T + \left[ b_1 \ \ldots \  b_n\right], \\ & b_i = \sum_{l=1}^\infty \frac{1}{l!} ([\Delta \vec{x}^T]^{\otimes l} \otimes I_n) \mbox{vec} (F_i^{(l+1)}(\theta)), \quad i=1, \ldots, n,\notag
\end{align} so that $b_i$ is a column vector and $F_q^{(l+1)}$ is as in \eqref{eq:Frecurrive}. Differentiating \eqref{eq:expZ} we may then write \eqref{eq:Zadjoint} as
\begin{align}\label{eq:Zbalance}
\sum_{k=0}^\infty \psi^k\left(\FD{Z^{(k)}}{t} + k\kappa Z^{(k)} \right) = -\left(J^T + \left[ b_1 \ \ldots \  b_n\right] \right) \sum_{k=0}^\infty \psi^k Z^{(k)}.
\end{align} Let $b_i^{(j)}$ denote the coefficient of $\psi^j$ in $b_i$. Then equating the coefficients of $\psi^k$ in \eqref{eq:Zbalance} we find that
\begin{align}\label{eq:Zk}
\FD{Z^{(k)}}{t} = -\left(J^T+ k\kappa Z^{(k)}\right)Z^{(k)} - \sum_{i=1}^n \sum_{j=0}^{k-1} {\rm e}_i^T Z^{(j)} b_i^{(k-j)}
\end{align} where ${\rm e}_i$ are the standard unit basis vectors. Using \eqref{eq:kroneckerDeltax} we can determine that
\begin{align}\label{eq:bvectors}
b_i^{(j)} = \sum_{l=1}^j \frac{1}{l!}\left( \left[ \sum_{\mathcal{A}^{l, j}} \bigotimes_{p=1}^l g^{(a_p)}\right]^T \otimes I_n\right)\mbox{vec}(F_i^{(l+1)}).
\end{align} Note that \eqref{eq:Zk}--\eqref{eq:bvectors} is in agreement with \cite{Wilson2020b} for a single isostable coordinate and also corrects for a missing term in the equation given there for $Z^{(2)}$.

Finally, in \cite{Wilson2020b} it is shown that the adjoint equation for $\mathcal{I}$ is
\begin{align}\label{eq:Iadjoint}
\FD{\mathcal{I}}{t} = -\left(\left.\PD{\vec{F}^T}{\vec{x}}\right|_{\vec{x}}- \kappa I_n\right) \mathcal{I}.
\end{align} Repeating the arguments for the expansion of $\mathcal{Z}$ above we see that
\begin{align}\label{eq:Ibalance} \begin{split}
    \sum_{k=0}^\infty \psi^k\left[\FD{I^{(k)}}{t} + k\kappa I^{(k)} \right] = -\left(J^T + \left[ b_1 \ \ldots \  b_n\right] - \kappa I_n\right) \\ \qquad \qquad \times \sum_{k=0}^\infty \psi^k I^{(k)}, \end{split}
\end{align} and then matching coefficients of $\psi^k$ yields \eqref{eq:IRCODE}.

\subsubsection{Normalisation of terms in PRC and IRC expansions}

Some of the equations \eqref{eq:PRCODE} and \eqref{eq:IRCODE} require normalisation to determine a unique periodic solution. If such a normalisation is required, it may also be computed from the asymptotic expansions and the definition of the phase and isostable coordinates. Following \cite{Wilson2020b} we find that for $Z^{(0)}$ the normalisation is $Z^{(0)}(\theta) \cdot \vec{F}(\vec{x}^\gamma(\theta))=\omega$ as expected. For $k\geq 1$ the normalisation for $Z^{(k)}$ (for a single isostable coordinate) can be found by equating powers of $\psi$ in the expansion of the equation $\omega = \dot{\theta} = \mathcal{Z}(\theta, \psi) \cdot \left(\vec{F}(\vec{x}^\gamma(\theta)) + \FD{\Delta \vec{x}}{t}\right)$ with the result that
\begin{align}
\label{eq:Znorm}\begin{split}
0= & Z^{(k)}(\theta)\cdot \vec{F}(\vec{x}^\gamma(\theta)) + \sum_{l=0}^{k-1} Z^{(l)}(\theta) \cdot J g^{(k-l)}(\theta) \\ & + \sum_{l=0}^{k-2} Z^{(l)}(\theta) \cdot \alpha^{(k-l)}(\theta).
\end{split}\end{align} The normalisations for $I^{(k)}$ can be found by equating powers of $\psi$ in the expansion of the equation $\kappa \psi = \dot{\psi}= \mathcal{I}(\theta, \psi) \cdot \left(\vec{F}(\vec{x}^\gamma(\theta)) + \FD{\Delta \vec{x}}{t}\right)$ resulting in the normalisation
\begin{align}
\begin{split} \label{eq:normI}
\delta_{1k} \kappa  &= I^{(k)}(\theta)\cdot \vec{F}(\vec{x}^\gamma(\theta)) + \sum_{l=0}^{k-1} I^{(l)}(\theta) \cdot J g^{(k-l)}(\theta) \\  & \ + \sum_{l=0}^{k-2} I^{(l)}(\theta) \cdot \alpha^{(k-l)}(\theta).\end{split}
\end{align}

A more useful normalisation can be determined for $I^{(0)}$ as follows. First observe that
\begin{align}
\left.\PD{\psi}{\vec{x}}\right|_{\vec{x}} = \left.\PD{\psi}{\vec{x}}\right|_{\vec{x}^\gamma} + \left.\PD{^2\psi}{\vec{x}^2}\right|_{\vec{x}^\gamma}\Delta \vec{x} + \cdots
\end{align} and hence expanding $\Delta \vec{x}$ as in \eqref{eq:expdeltax} and comparing powers of $\psi$ with \eqref{eq:expI} we see that $I^{(0)} = \left.\PD{\psi}{\vec{x}}\right|_{\vec{x}^\gamma}$ and $I^{(1)} = \left.\PD{^2\psi}{\vec{x}^2}\right|_{\vec{x}^\gamma} g^{(1)}$. Then $\FD{I^{(0)}}{t} = \left.\PD{^2\psi}{\vec{x}^2}\right|_{\vec{x}^\gamma(\theta)}^T \vec{F}(\vec{x}^\gamma(\theta))$. Therefore \eqref{eq:normI} becomes
\begin{align*}
\kappa &=  \left( \left.\PD{^2\psi}{\vec{x}^2}\right|_{\vec{x}^\gamma} g^{(1)}\right)\cdot \vec{F}(\vec{x}^\gamma(\theta)) + I^{(0)}(\theta)\cdot  J g^{(1)}(\theta)\\
& = g^{(1)}(\theta) \cdot  \left.\PD{^2\psi}{\vec{x}^2}\right|_{\vec{x}^\gamma}^T \vec{F}(\vec{x}^\gamma(\theta)) + g^{(1)}(\theta) \cdot J^T I^{(0)}(\theta)\\
& = g^{(1)}(\theta) \cdot  \left(\FD{I^{(0)}}{t} + J^T I^{(0)}(\theta)\right) = \kappa g^{(1)}(\theta) \cdot I^{(0)}(\theta),
\end{align*} where \eqref{eq:IRCODE} is used in the final step.
Thus the required normalisation is satisfied if $I^{(0)}(0)\cdot g^{(1)}(0) =1$.

\section{Stability of phase-locked states}\label{sec:AppendixB}

 \subsection{General phase-locked states}
For a phase locked state $\Phi= (\phi_1, \ldots, \phi_N)$ the linearisation of \eqref{eqs:piaveraged} has $2N \times 2N$ block matrix Jacobian
\begin{align*}\mathcal{J} =
\left[\begin{array}{@{}c|c@{}}
  \mathcal{H}^{(1)}(\Phi) & \mathcal{H}^{(2)}(\Phi) \\
\hline
  \mathcal{H}^{(3)}(\Phi) & \mathcal{H}^{(4)}(\Phi)
\end{array}\right]
\end{align*} where \begin{subequations}\label{eqs:jacobian}
    \begin{align}
\begin{split}\label{eq:Hij1}
\mathcal{H}_{ij}^{(1)}(\Phi) = & \ \epsilon w_{ij}\Xi^{(\theta)}_{ij} -\epsilon\delta_{ij}\sum_{k=1}^N w_{ik}\Xi^{(\theta)}_{ik},\end{split}\\ \label{eq:Hij2}
\mathcal{H}_{ij}^{(2)}(\Phi) = &\  \epsilon w_{ij} H_3(\phi_j-\phi_i) + \epsilon\delta_{ij}\sum_{k=1}^N w_{ik} H_2(\phi_k - \phi_i),\\ \label{eq:Hij3}
\begin{split}
\mathcal{H}_{ij}^{(3)}(\Phi) = &\  \epsilon w_{ij}\Xi^{(\psi)}_{ij} -\epsilon\delta_{ij}\sum_{k=1}^N w_{ik}\Xi^{(\psi)}_{ik} ,\end{split}
\\ \label{eq:Hij4}
\begin{split}
\mathcal{H}_{ij}^{(4)}(\Phi) = &\  \epsilon w_{ij} H_6(\phi_j-\phi_i) + \delta_{ij}\left(\kappa+\varepsilon\sum_{k=1}^N w_{ik} H_5(\phi_k - \phi_i)\right)
\end{split}
\end{align}\end{subequations} for 
\begin{align*}
   \Xi^{(\theta)}_{ij}  & = H'_{1}(\phi_j - \phi_i)+\Psi_iH'_{2}(\phi_j - \phi_i)+\Psi_jH'_{3}(\phi_j - \phi_i),\\
   \Xi^{(\psi)}_{ij}  & = H'_{4}(\phi_j - \phi_i)+\Psi_iH'_{5}(\phi_j - \phi_i)+\Psi_jH'_{6}(\phi_j - \phi_i).
\end{align*}

\subsection{Two-cluster states}
For two cluster states where the two clusters have distinct orbits, the linearisation about $\Phi = (0, \ldots, 0, \chi, \ldots, \chi)$ has $2N \times 2N$ block matrix Jacobian \eqref{eq:pljac} where each of the blocks $\mathcal{H}^{(q)}$ are themselves block matrices of the form \eqref{eq:2clusterblocks}. The coefficients in the definition of the blocks \eqref{eq:blockAA}--\eqref{eq:blockBB} can be rearranged into arrays and are found to be given by 
\begin{widetext}
\begin{subequations}
    \begin{align} \begin{split}
\mathcal{A} &= \begin{bmatrix} a_1 & a_2\\ a_3 &a_4\end{bmatrix} = \frac{\epsilon}{N} \begin{bmatrix} H_1'(0) + {\Psi}_A(H_2'(0) + H_3'(0)) & H_3(0) \\H_4'(0) + {\Psi}_A(H_5'(0) + H_6'(0)) & H_6(0)   \end{bmatrix}, \end{split}\\ \begin{split}
\mathcal{B}&=  \begin{bmatrix} b_1 & b_2\\ b_3 &b_4\end{bmatrix} = \frac{\epsilon}{N} \begin{bmatrix} H_1'(\chi) + {\Psi}_A H_2'(\chi) + {\Psi}_B H_3'(\chi) & H_3(\chi) \\H_4'(\chi) + {\Psi}_A H_5'(\chi) + {\Psi}_B H_6'(\chi) &  H_6(\chi)   \end{bmatrix}, \end{split}\\ \begin{split}
\mathcal{C}&=  \begin{bmatrix} c_1 & c_2\\ c_3 &c_4\end{bmatrix} = \frac{\epsilon}{N} \begin{bmatrix} H_1'(-\chi) + {\Psi}_B H_2'(-\chi) + {\Psi}_A H_3'(-\chi) & H_3(-\chi) \\H_4'(-\chi) + {\Psi}_B H_5'(-\chi) + {\Psi}_A H_6'(-\chi) &  H_6(-\chi)   \end{bmatrix}, \end{split}\\ \begin{split}
\mathcal{D} &= \begin{bmatrix} d_1 & d_2\\ d_3 &d_4\end{bmatrix} = \frac{\epsilon}{N} \begin{bmatrix} H_1'(0) + {\Psi}_B(H_2'(0) + H_3'(0)) & H_3(0) \\H_4'(0) + {\Psi}_B(H_5'(0) + H_6'(0)) & H_6(0)   \end{bmatrix}, \end{split}\\ \begin{split}
\mathcal{S^A} &= \begin{bmatrix} s^A_1 & s^A_2\\ s^A_3 &s^A_4\end{bmatrix} =  \begin{bmatrix} -N_A a_1 -N_B b_1 & \epsilon (N_AH_2(0) + N_B H_2(\chi))/N \\-N_A c_1 - N_Bd_1  & \kappa + \epsilon ( N_A H_5(0) +  N_B H_5(\chi))/N\end{bmatrix}, \end{split}\\ \begin{split}
\mathcal{S^B} &= \begin{bmatrix} s^B_1 & s^B_2\\ s^B_3 &s^B_4\end{bmatrix} =  \begin{bmatrix} -N_A a_3 -N_B b_3 & \epsilon ( N_A H_2(-\chi) + N_B H_2(0))/N \\-N_A c_3 - N_B d_3  & \kappa +\epsilon ( N_A H_5(-\chi) + N_B H_5(0))/N \end{bmatrix}. \end{split}
\end{align}
\end{subequations}

We can then rewrite
\begin{align} 
|\mathcal{J}- \mu I_{2N}| = \begin{vmatrix}
\mathcal{A} \otimes \mathbbm{1}_{N_A \times N_A} + (\mathcal{S}^A-\mu I_2) \otimes I_{N_A} & \mathcal{B} \otimes \mathbbm{1}_{N_A \times N_B} \\ \mathcal{C} \otimes \mathbbm{1}_{N_B \times N_A} & \mathcal{D} \otimes \mathbbm{1}_{N_B \times N_B} + (\mathcal{S}^B-\mu I_2) \otimes I_{N_B}
\end{vmatrix}.
\end{align}
\end{widetext}

A generalisation of the matrix determinant lemma \cite{Ding2007} says that for invertible $n \times n$ matrix $P$ and invertible $m \times m$ matrix $Q$, if $U$ and $V$ are $n \times m$ matrices then 
\begin{align} \label{eq:mdl}
\det(P + UQV^T) = \det(Q^{-1} + V^T P^{-1} U) \det(Q) \det(P).
\end{align}
Taking $P=(\mathcal{S}^A-\mu I_2)\otimes I_{N_A}$, $Q= (\mathcal{S}^A-\mu I_2)^{-1}$ and letting $U=(\mathcal{S}^A-\mu I_2) \otimes \mathbbm{1}_{N_A \times 1}$ and $V= \mathcal{A}^T \otimes \mathbbm{1}_{N_A\times 1}$ we observe that since $UQV^T= \mathcal{A} \otimes \mathbbm{1}_{N_A \times N_A}$ and $V^TP^{-1}U = N_A \mathcal{A}$, 
\begin{align} \begin{split}
& \det(\mathcal{A} \otimes \mathbbm{1}_{N_A \times N_A} + (\mathcal{S}^A-\mu I_2) \otimes I_{N_A}) \\ & \ \   = \det(\mathcal{S}^A-\mu I_2)^{N_A-1} \det(\mathcal{S}^A + N_A \mathcal{A} -\mu I_2).\end{split}
\end{align} This reveals that eigenvalues corresponding to intracluster perturbations are given by the eigenvalues of $\mathcal{S}^A$ each with multiplicity $N_A-1$ and the eigenvalues  of $\mathcal{S}^B$ each with multiplicity $N_B-1$. The remaining four eigenvalues can be found by considering intercluster perturbations $(\phi_A, \psi_A, \phi_B, \psi_B)= (\overline{\phi}_A, {\Psi}_A, \overline{\phi}_B, {\Psi}_B) + (\widetilde{\phi}_A, \widetilde{\psi}_A, \widetilde{\phi}_B, \widetilde{\psi}_B)$ where $\overline{\phi}_A$ and $\overline{\phi}_B$ are the phases of the two clusters in the phase locked state. Then 
\begin{align}
\FD{}{t} \begin{bmatrix}
\widetilde{\phi}_A\\ \widetilde{\psi}_A\\ \widetilde{\phi}_B\\ \widetilde{\psi}_B
\end{bmatrix} =\begin{bmatrix}\mathcal{S}^A + N_A \mathcal{A} & N_B \mathcal{B}\\ N_A \mathcal{C} & \mathcal{S}^B + N_B \mathcal{D} 
\end{bmatrix} \begin{bmatrix}
\widetilde{\phi}_A\\ \widetilde{\psi}_A\\ \widetilde{\phi}_B\\ \widetilde{\psi}_B
\end{bmatrix}.\label{eq:intercluter}
\end{align} Therefore the intercluster eigenvalues are the eigenvalues of the Jacobian $\mathcal{J}_M$ in \eqref{eq:intercluter} which is also given in \eqref{eq:Jtwoclust}. One of these is zero corresponding to pure rotations, observing that the first and third columns of $\mathcal{J}_M$ are identical. The others correspond to changes in the phase difference between clusters and changes in the isostable values for the two synchronised clusters. In conclusion, 
\begin{align}\begin{split}
\det(\mathcal{J}- \mu I_{2N}) =&  \det( \mathcal{S}^A - \mu I_2)^{N_A -1} \det( \mathcal{S}^B - \mu I_2)^{N_B -1} \\ &\qquad \times  \det(\mathcal{J}_M - \mu I_4). \end{split}
\end{align}  

\section{Higher order phase equation derivation}\label{sec:AppendixC}

Here we follow and correct Park and Wilson \cite{Park2021} to derive the higher order phase equation \eqref{eq:higerorderphase} explicitly including terms up to cubic order in the interaction strength $\epsilon$. Park and Wilson \cite{Park2021} take the isostable coordinate as an $O(\epsilon)$ term which may be expressed as $\psi(t)= \epsilon p^{(1)}(t) + \epsilon^2 p^{(2)}(t) + \cdots$ where $p^{(k)}(t)$ are $O(1)$. Using this additional expansion in \eqref{eq:phaseiso2} and matching terms at different powers of $\epsilon$ results in a hierarchy of linear first order differential equations for the $p^{(k)}(t)$ \cite{Park2021}. The equations at $O(\epsilon)$ and $O(\epsilon^2)$ are
\begin{align}
\FD{p_i^{(1)}}{t} &= \kappa p_i^{(1)} + \sum_{j=1}^N w_{ij} h_4(\theta_i, \theta_j)\\
\FD{p_i^{(2)}}{t} &= \kappa p_i^{(2)} + \sum_{j=1}^N w_{ij} \left[ p_i^{(1)} h_5(\theta_i, \theta_j) + p_j^{(1)} h_6(\theta_i, \theta_j)\right],
\end{align} where $h_4, h_5, h_6$ are given by \eqref{eq:h4}--\eqref{eq:h6}. The forcing terms of the equation at a given order depend only on the solutions of lower order equations and hence each can be solved explicitly in turn using the integrating factor $e^{-\kappa t}$. Following \cite{Park2021, Wilson2019c} (and also correcting notational errors in \cite{Park2021}) we observe that
\begin{align}
p_i^{(1)}(t) &= \sum_{k=1}^N w_{ik} q^{(1)}(\theta_i, \theta_k)\\
p_i^{(2)}(t) &= \sum_{k=1}^N\sum_{l=1}^N \left[w_{ik}w_{il} q^{(2)}(\theta_i, \theta_k, \theta_l) + w_{ik}w_{kl} q^{(3)}(\theta_i, \theta_k, \theta_l)\right] ,
\end{align} where
\begin{subequations}
\begin{align}
q^{(1)}(\theta_i, \theta_k) &= \int_0^\infty e^{\kappa s } h_4(\theta_i-\omega s, \theta_k-\omega s)\ {\rm d} s,\\
\begin{split}
q^{(2)}(\theta_i, \theta_k, \theta_l) &= \int_0^\infty e^{\kappa s } q^{(1)}(\theta_i-\omega s, \theta_l-\omega s) \\ & \qquad \times h_5(\theta_i-\omega s, \theta_k-\omega s)\ {\rm d} s, \end{split}\\ \begin{split}
q^{(3)}(\theta_i, \theta_k, \theta_l) &= \int_0^\infty e^{\kappa s } q^{(1)}(\theta_k- \omega s, \theta_l-\omega s) \\ & \qquad \times h_6(\theta_i-\omega s, \theta_k-\omega s)\ {\rm d} s.\end{split}
\end{align}
\end{subequations} Higher order terms $p^{(k)}$, $k>2$ may also be determined (see \cite{Park2021}) but we do not require these here. Substituting all expansions into \eqref{eq:phaseiso} we obtain the higher order phase equation
\begin{widetext}
    \begin{align}\label{eq:higerorderphase}
\begin{split}
\FD{\theta_i}{t} = \omega \ +  \ & \epsilon \sum_{j=1}^N w_{ij} \overline{h}_1(\theta_i, \theta_j) +  \epsilon^2 \sum_{j,k=1}^N \left[ w_{ij}w_{ik} \overline{h}_2(\theta_i, \theta_j, \theta_k) + w_{ij}w_{jk} \overline{h}_3(\theta_i, \theta_j, \theta_k) \right] \\  + \  & \epsilon^3 \sum_{j,k, l=1}^N \Bigl[ w_{ij}w_{ik}w_{il} \overline{h}_4(\theta_i, \theta_j, \theta_k, \theta_l)  + w_{ij}w_{ik}w_{kl}  \overline{h}_5(\theta_i, \theta_j, \theta_k, \theta_l)  + w_{ij}w_{jk}w_{jl}  \overline{h}_6(\theta_i, \theta_j, \theta_k, \theta_l) \\ & \qquad \qquad \qquad \qquad + w_{ij}w_{jk}w_{kl}  \overline{h}_7(\theta_i, \theta_j, \theta_k, \theta_l)  + w_{ij}w_{ik}w_{jl}  \overline{h}_8(\theta_i, \theta_j, \theta_k, \theta_l)\Bigr],
\end{split}
\end{align}\end{widetext} where
\begin{subequations}
    \begin{align}\label{eq:hbar1}
\overline{h}_1(\theta_i, \theta_j) &= h_1(\theta_i, \theta_j),\\
\overline{h}_2(\theta_i, \theta_j, \theta_k) &= q^{(1)}(\theta_i, \theta_k) h_2(\theta_i, \theta_j),\\
\overline{h}_3(\theta_i, \theta_j, \theta_k) &= q^{(1)}(\theta_j, \theta_k) h_3(\theta_i, \theta_j),\\
\begin{split}
\overline{h}_4(\theta_i, \theta_j, \theta_k, \theta_l) &= q^{(2)}(\theta_i, \theta_k, \theta_l)h_2(\theta_i, \theta_j) \\ & \quad + q^{(1)}(\theta_i, \theta_k) q^{(1)}(\theta_i, \theta_l) h_7(\theta_i, \theta_j),\end{split}\\
\overline{h}_5(\theta_i, \theta_j, \theta_k, \theta_l) &= q^{(3)}(\theta_i, \theta_k, \theta_l)h_2(\theta_i, \theta_j), \\
\begin{split}
\overline{h}_6(\theta_i, \theta_j, \theta_k, \theta_l) &= q^{(2)}(\theta_j, \theta_k, \theta_l)h_3(\theta_i, \theta_j) \\ & \quad +  q^{(1)}(\theta_j, \theta_k) q^{(1)}(\theta_j, \theta_l) h_8(\theta_i, \theta_j), \end{split}\\
\overline{h}_7(\theta_i, \theta_j, \theta_k, \theta_l) &=q^{(3)}(\theta_j, \theta_k, \theta_l)h_3(\theta_i, \theta_j), \\
\overline{h}_8(\theta_i, \theta_j, \theta_k, \theta_l) &=  q^{(1)}(\theta_i, \theta_k) q^{(1)}(\theta_j, \theta_l) h_9(\theta_i, \theta_j),\label{eq:hbar8}
\end{align} \end{subequations}
for $h_1, \ldots, h_6$ as defined in \eqref{eqs:hs}, and for
\begin{subequations}
    \begin{align} \begin{split}
h_7(\theta_i, \theta_j)& = Z^{(0)}(\theta_i) \cdot K_1(\theta_i, \theta_j) + Z^{(1)}(\theta_i) \cdot J_1 g^{(1)}(\theta_i) \\ & \quad + Z^{(2)}(\theta_i) \cdot \vec{G}(\vec{x}^\gamma(\theta_i), \vec{x}^\gamma(\theta_j)), \end{split}\\
h_8(\theta_i, \theta_j)& = Z^{(0)}(\theta_i) \cdot K_2(\theta_i, \theta_j), \\
h_9(\theta_i, \theta_j)& =Z^{(0)}(\theta_i) \cdot L(\theta_i, \theta_j) + Z^{(1)}(\theta_i) \cdot J_2 g^{(1)}(\theta_j).
\end{align} \end{subequations}
Moving to a rotating reference frame, applying first-order averaging and returning to the original variables we can obtain the autonomous approximation to the higher order phase reduced equation \eqref{eq:avhigherorderphase} where 
\begin{subequations}
\begin{align}
\overline{H}_1(\chi) &= \frac{1}{2\pi}\int_0^{2\pi} \overline{h}_1(u, u+\chi) \,{\rm d}u,\\
\overline{H}_m(\chi, \eta) &= \frac{1}{2\pi}\int_0^{2\pi} \overline{h}_m(u, u+\chi, u+\eta) \,{\rm d}u, \quad m=2,3,\\
\begin{split}
\overline{H}_m(\chi, \eta, \xi) &= \frac{1}{2\pi}\int_0^{2\pi} \overline{h}_m(u, u+\chi, u+\eta, u+ \xi) \,{\rm d}u, \\ & \qquad \qquad \qquad \qquad m=4,\ldots, 8.\end{split}
\end{align}
\end{subequations}

\section{Splay state linear stability boundary for $N\geq 3$}\label{sec:AppendixD}

  In the large $N$ limit, from \eqref{eqs:lambda} we calculate
\begin{align}\label{eq:lambda1MFCGL}
\lambda_q^{(1)}  =& \begin{cases} -\epsilon i \left(i(1+c_1c_2) + c_2-c_1 \right)/2 & \text{for } q=1, \\  0 & \text{otherwise, }\end{cases}\\
\lambda_q^{(2)}  = & \begin{cases} A(i-c_1)(1+c_2^2)\epsilon/2 & \text{for } q=1, \\  0 & \text{otherwise, }\end{cases}\\
\lambda_q^{(3)}  = & \begin{cases} \frac{\epsilon i (ic_1+1)}{2A(\epsilon-1)} & \text{for } q=1, \\  0 & \text{otherwise, }\end{cases}\\
\lambda_q^{(4)}  = & \begin{cases} \kappa + \epsilon\left((i-c_1)c_2 + ic_1 +5\right)/2 & \text{for } q=1, \\  2(\epsilon-1) & \text{otherwise. }\end{cases} \label{eq:lambda4MFCGL}
\end{align} Therefore using \eqref{eq:splayeigs}, the eigenvalues of the Jacobian for the splay state are $0$ and $2(\epsilon -1)$, each of multiplicity $N-1$, and the eigenvalues of $\Lambda_1$. It can be shown that there is a pair of eigenvalues $\mu=\pm i \sigma$ with
\begin{align}
\sigma = \frac{\epsilon\left(c_2^2 c_1 - 2c_2 + 3c_1)\epsilon^2 + 2(c_2-c1)(2\epsilon -1)\right)}{2(\epsilon-1)(3\epsilon -2)},
\end{align} when $\epsilon= \epsilon_{0, PI}$ satisfies \eqref{eq:splaystabboundary}. Equation \eqref{eq:splaystabboundary} defines a stability boundary for the splay state in the large $N$ limit. It can similarly be shown that for finite $N\geq 3$ the eigenvalues of the Jacobian are $0$, $2(\epsilon-1)$ each of multiplicity $N-2$ and the remaining eigenvalues are the eigenvalues of $\Lambda_1$ and their complex conjugates. This gives the identical equation \eqref{eq:splaystabboundary} for the stability boundary.

\section{The Morris-Lecar model} \label{sec:AppendixE}

The Morris-Lecar model is a planar model of neural activity \cite{Morris1981} with two ionic currents; an outward non-inactivating potassium current and an inward, non-inactivating calcium current. Assuming that the calcium dynamics operate on a much faster time scale than the potassium, the dynamics of an isolated node are given by
\begin{subequations}\label{eqs:ML}
    \begin{align}
\begin{split}
C_m \FD{v_i}{t} &= I_b - g_L(v_i - E_L)-g_Kw_i(v_i - E_K)\\ & \quad -g_{Ca}m_{\infty}(v_i)(v_i-E_{Ca}) \label{eq:MLv} \end{split}\\
\FD{w_i}{t} &= \phi (w_{\infty}(v_i)-w_i)\lambda(v_i), \label{eq:MLw}
\end{align}
\end{subequations}
where
\begin{align*}
m_{\infty}(v) &= 0.5\left(1+\tanh\left(\left(v -V_1\right)/V_2\right)\right) \\
w_{\infty}(v) &= 0.5\left(1+\tanh\left(\left(v-V_3\right)/V_4\right)\right)\\
\lambda(v) &= \cosh \left(\left(v-V_3\right)/\left(2V_4\right)\right).
\end{align*} and $V_1, \ldots, V_4$ and $\phi$ are constants. Here $w_i$ represents the fraction of K$^+$ channels open at node $i$ and the Ca$^{2+}$ channels respond to $v_i$ so rapidly that we assume instant activation. Here $g_L$ is the leakage conductance,  $g_K$ and $g_{Ca}$ are the potassium and calcium conductances, $E_L$, $E_K$, $E_{Ca}$ are corresponding reversal potentials, $m_\infty(v)$, $w_\infty(v)$ are voltage dependent gating functions, $\lambda(v)$ is a voltage-dependent rate, $I_b$ is the applied current and $C_m$ is the cell capacitance. Unless otherwise stated, the default model parameters used throughout are as in \cite{Han1995}: $\phi = 1.15$, $g_{Ca} = 1$, $g_K = 2$, $g_L = 0.5$, $E_{Ca} = 1$, $E_K = -0.7$, $E_L = -0.5$, $V_1 = -0.01$, $V_2 = 0.15$, $V_3 = 0.1$, $V_4 = 0.145$ and $C_m = 1$. For these parameter values a limit cycle of the system arises at $I_b = 0.0730$ through a homoclinic connection. We choose $I_b= 0.075$ to place the system near to the homoclinic bifurcation. For these parameter values we find that the periodic orbit has Floquet exponent $\kappa = -0.4094$ and period $T=8.1654$.

\bibliography{PINetworks}

\end{document}